\documentclass[mathpazo]{cicp}

\usepackage{times,amsmath,amsbsy,amssymb,amscd,mathrsfs}
\usepackage{mathtools}
\usepackage{algorithm2e} 

\usepackage{multicol,multirow}
\usepackage{booktabs}

\usepackage[labelformat=empty]{subfig}
\usepackage{tikz,tikz-cd}
\usepackage{pgfplots}
\pgfplotsset{compat=newest}
\pgfplotsset{plot coordinates/math parser=false}
\newlength\fwidth

\definecolor{myBlue}{rgb}{0.0,0.0,0.55}
\usepackage[pdftex,colorlinks=true,citecolor=myBlue,linkcolor=myBlue]{hyperref}

\usepackage{enumerate,xspace}

\newcommand{\bs}{\boldsymbol}

\newcommand{\vertiii}[1]{{\left\vert\kern-0.25ex\left\vert\kern-0.25ex\left\vert #1 \right\vert\kern-0.25ex\right\vert\kern-0.25ex\right\vert}}

\usepackage{algorithm}
\usepackage{algpseudocode}
\usepackage{stmaryrd}
\usepackage{scalefnt}
\usepackage{graphicx}
\usepackage{caption}
\usepackage{listings}
\usepackage[T1]{fontenc}
\usepackage{siunitx}
\usepackage{float}

\definecolor{codegray}{rgb}{0.95,0.95,0.95}
\definecolor{codeframe}{rgb}{0.5,0.5,0.5}
\definecolor{codecomment}{rgb}{0.4,0.4,0.4}
\definecolor{codegreen}{rgb}{0,0.6,0}
\definecolor{codeblue}{rgb}{0,0,0.8}
\definecolor{codepurple}{rgb}{0.6,0,0.6}

\begin{document}
\title[Efficient Multi-Backend Implementation]{SOPTX: A High-Performance Multi-Backend Framework for Topology Optimization}


\author[He L, Wei H Y, Tian T]{Liang He \affil{1}, 
Huayi Wei\affil{2}\comma\corrauth, Tian Tian\affil{3}}
\address{
\affilnum{1}\ 
School of Mathematics and Computational Science, Xiangtan University, Xiangtan 411105, China\\
\affilnum{2}\ 
School of Mathematics and Computational Science, Xiangtan University; 
National Center of  Applied Mathematics in Hunan, Hunan Key Laboratory for 
Computation and Simulation in Science and Engineering, Xiangtan 411105, China\\
\affilnum{3}\ 
School of Mathematics and Computational Science, Xiangtan University, Xiangtan 411105, China
}

\emails{
	{\tt lianghe@smail.xtu.edu.cn} (L.~He), 
	{\tt weihuayi@xtu.edu.cn} (H.~Wei),
	{\tt tiantian@smail.xtu.edu.cn} (T.~Tian)
	}


\begin{abstract}
In recent years, topology optimization (TO) has gained widespread attention as a powerful structural design method. However, its application remains challenging due to the deep expertise and extensive development effort required. Traditional TO methods, tightly coupled with computational mechanics like finite element method (FEM), result in intrusive algorithms demanding a comprehensive system understanding. This paper presents SOPTX, a TO package based on FEALPy, which implements a modular architecture that decouples analysis from optimization, supports multiple computational backends (NumPy, PyTorch, JAX), and achieves a non-intrusive design paradigm. Core innovations include: (1)
cross-platform design that supports multiple computational backends, enabling efficient algorithm execution on central processing units (CPUs) and flexible acceleration using graphics processing units (GPUs), while leveraging automatic differentiation (AD) technology for efficient sensitivity computation of objective and constraint functions; (2) fast matrix assembly techniques that overcome the performance bottlenecks of traditional numerical integration methods, significantly accelerating finite element computations and enhancing overall efficiency; (3) a modular framework supporting TO problems for arbitrary dimensions and meshes, allowing flexible configuration and extensibility of optimization workflows through a rich library of composable components. Using the density-based method for the classic compliance minimization problem with volume constraints as an example, numerical experiments demonstrate SOPTX's high efficiency in computational speed and memory usage, while showcasing its strong potential for research and engineering applications.
\end{abstract}


\maketitle

\section{Introduction}
Topology optimization (TO) is a class of structural optimization techniques aimed at improving structural performance by optimizing material distribution within a given design domain. Widely applied in aerospace, automotive, and civil engineering, TO addresses critical design challenges by efficiently utilizing materials. Among various approaches, density-based methods are particularly popular due to their intuitive concept and practicality. The Solid Isotropic Material with Penalization (SIMP) method, introduced by Bendsøe et al.~\cite{Bendsøe2004}, is the most widely adopted density-based approach, promoting binary (0 - 1) solutions by penalizing intermediate densities and integrating effectively with finite element method (FEM).

However, TO inherently involves large-scale computations due to the tight coupling between structural analysis and optimization. Each iteration requires solving boundary value problems and performing sensitivity analysis for gradient-based optimization algorithms. For large-scale problems, repeatedly solving large linear systems and evaluating sensitivities imposes significant demands on computational resources. Therefore, improving computational efficiency and scalability without sacrificing accuracy remains a major challenge in TO research and applications.

To lower the entry barrier and promote widespread adoption of TO, researchers have published numerous educational studies and literature. A pioneering example is the 99-line MATLAB code by Sigmund~\cite{sigmund200199}, which demonstrated the fundamentals of a two-dimensional SIMP algorithm in a concise and self-contained manner, profoundly impacting both TO education and research. Subsequently, improved versions of this educational code have emerged, including the more efficient 88-line version proposed by Andreassen et al.~\cite{andreassen2011efficient} and the extended 3D implementation by Liu and Tovar~\cite{liu2014efficient}. These codes convey practical TO knowledge simply and provide self-contained examples of basic numerical algorithms.

Meanwhile, efforts in open-source software development have also advanced TO practices. Chung et al.~\cite{chung2019topology} leveraged OpenMDAO~\cite{gray2019openmdao}, decomposing TO into modular components with automated derivative assembly, enhancing flexibility and extensibility. Gupta et al.~\cite{gupta202055} developed a parallel-enabled implementation using the FEniCS finite element framework~\cite{alnaes2015fenics}, addressing large-scale optimization. Recently, Ferro and Pavanello~\cite{ferro2023simple} provided an efficient 51-line TO implementation integrating FEniCS, Dolfin Adjoint, and Interior Point Optimizer, simplifying the development process. These projects provided standardized interfaces for convenient integration with existing finite element analysis (FEA) tools, reducing implementation complexity. However, despite improvements, existing open-source TO packages still face limitations in terms of functionality extensibility and flexibility, particularly for complex engineering applications.

To accelerate the development of TO, automating sensitivity analysis has become a critical step. This involves automatically computing derivatives of objectives, constraints, material models, projections, filters, and other components with respect to the design variables. Currently, the common practice involves manually calculating sensitivities, which can be tedious and error-prone despite their theoretical simplicity, often becoming a bottleneck in the development of new TO modules and exploratory research. Automatic differentiation (AD) provides an efficient and accurate approach for evaluating derivatives of numerical functions ~\cite{griewank2008evaluating}. By decomposing complex functions into a series of elementary operations (such as addition and multiplication), AD accurately computes derivatives of arbitrary differentiable functions. In TO, the Jacobian matrices represent sensitivities of objective functions and constraints with respect to design variables, and software can automate this process, relieving developers from manually deriving and implementing sensitivity calculations. With its capability of easily obtaining accurate derivative information, AD offers significant advantages in design optimization, particularly for highly nonlinear optimization problems.

In recent years, the adoption of AD in TO has gradually increased. For instance, Nørgaard et al.~\cite{norgaard2017applications} employed the AD tools CoDiPack and Tapenade to achieve automatic sensitivity analysis in unsteady flow TO, significantly enhancing computational efficiency. Building upon this, Chandrasekhar~\cite{chandrasekhar2021auto} utilized the high-performance Python library JAX~\cite{jax2018github} to apply AD techniques in density-based TO, achieving efficient solutions to classical TO problems, exemplified by compliance minimization.

Although existing TO programs, ranging from educational codes to open-source implementations, have significantly advanced the field, they often employ a procedural programming paradigm that tightly couples analysis modules, such as FEA, with optimization modules. This tight coupling restricts software extensibility, necessitating invasive changes for new features, and diminishes reusability, complicating integration into broader frameworks like multidisciplinary design optimization (MDO) in aerospace. Consequently, developing a decoupled architecture that enhances both extensibility and reusability without sacrificing accuracy remains a critical challenge in TO.

To address the aforementioned challenges, this paper proposes SOPTX, a high-performance TO framework built upon the FEALPy platform ~\cite{fealpy}. FEALPy is an open-source intelligent CAX computing engine providing an efficient, flexible, and extensible platform for numerical simulation. Its selection as the underlying platform is supported by key advantages:
\vspace{-1ex}
\begin{enumerate}
	\item FEALPy supports multiple computational backends, including NumPy, PyTorch, and JAX, enabling efficient execution across a wide range of hardware architectures such as central processing units (CPUs) and graphics processing units (GPUs).
	\item FEALPy supports a broad range of numerical schemes, including FEM of arbitrary order and dimension, as well as alternative methods such as the virtual element method (VEM) and the finite difference method (FDM).
	\item FEALPy adopts highly efficient vectorized operations that fully exploit the computational power of modern processors.
	\item FEALPy features a clean and extensible API design, which aligns well with our goal of developing a modular and reusable TO framework.
\end{enumerate}
\vspace{-1ex}
Consequently, FEALPy provides a robust and versatile foundation for developing modular and extensible TO frameworks.

Leveraging the strengths of FEALPy, the SOPTX framework achieves a highly modular design. SOPTX adopts a component-based philosophy commonly used in MDO, decomposing complex systems into independent and reusable modules to enhance flexibility and maintainability. Specifically, SOPTX introduces a clear separation between analysis and optimization: numerical solvers, AD tools, and optimization algorithms are designed as independent and interchangeable modules. This design allows users to flexibly select or replace individual modules based on specific needs without requiring substantial changes to the overall code structure. Through this modular architecture, SOPTX not only inherits the efficiency and flexibility of FEALPy but also significantly enhances the framework's extensibility and reusability, providing robust support for complex engineering applications.

Building on this modular foundation, SOPTX seamlessly integrates AD into the TO workflow. It supports multiple computational backends and can automatically switch between them based on hardware availability and performance requirements. Additionally, SOPTX mitigates the computational bottlenecks associated with traditional numerical integration through a fast matrix assembly technique. This technique separates element-dependent and element-independent components, avoiding redundant computations of invariant data during iterative procedures. Furthermore, symbolic integration via SymPy~\cite{meurer2017sympy} replaces conventional numerical quadrature, reducing computational cost while improving accuracy. SOPTX employs a non-intrusive design, allowing users to easily replace or extend filters, constraints, and objectives without modifying the core framework. Overall, SOPTX is designed as a low-barrier, high-performance platform that, through the combination of modularity and AD, empowers users to focus on algorithmic innovation and problem modeling rather than implementation details. This makes it well-suited for education, research, and complex multiphysics-coupled scenarios in TO and MDO.

The remainder of this paper is organized as follows. Section 2 introduces the mathematical formulation of TO, covering density-based methods, compliance minimization, optimization algorithms, and filtering techniques. Section 3 presents the design philosophy of SOPTX built on FEALPy, with emphasis on its modular architecture and backend switching for high-performance, flexible implementation. Section 4 demonstrates SOPTX installation and usage via a classical 2D cantilever beam problem. Section 5 showcases its effectiveness through a series of 2D and 3D examples, including comparative studies with different filters and optimization algorithms, and performance analyses of fast matrix assembly, AD, and backend switching. Section 6 concludes the paper and outlines future research directions.

\section{Density-based Topology Optimization: Formulation, Algorithms, and Regularization}\label{sec:math}
In this section, we introduce the density-based method in topology optimization (TO), explaining how it optimizes material distribution using a density function and interpolation models. Using the compliance minimization problem as an example, we present its continuous and discrete formulations. We then describe two key optimization algorithms in TO: the Optimality Criteria (OC) method and the Method of Moving Asymptotes (MMA), including their implementations and features. Finally, we discuss key filtering techniques in TO, such as sensitivity, density, and Heaviside projection filters, and their role in improving numerical stability and physical feasibility.

\subsection{Density-based Method}\label{sec:DbM}
The core objective of TO is to determine the optimal material distribution within a given design domain to achieve predefined performance targets. Density-based methods are among the most widely used strategies in the field of TO. These methods parameterize the structure by defining a material density function $\rho(x)$, where $\rho(x) = 1$ represents solid material, $\rho(x) = 0$ represents voids, and intermediate values $0 < \rho(x) < 1$ are used to handle the non-convexity of discrete problems.

In a discretized design domain, the mechanical properties (e.g., stiffness) of material elements transition smoothly between solid and void states via interpolation models~\cite{bendsoe1989optimal}. A widely used model is the Solid Isotropic Material with Penalization (SIMP) approach~\cite{zhou1991coc}, establishes a power-law relationship between the material density $\rho$ and Young's modulus $E$:
\begin{equation*}
	E(\rho) = \rho^p E_0,\quad \rho \in [0,1],
\end{equation*}
where $E_0$ denotes the Young's modulus of solid material, and $p>1$ is the penalization factor. By increasing the relative stiffness cost of intermediate densities, this power-law relationship encourages the optimization results toward a clear 0-1 distribution.

To avoid singularities in the stiffness matrix, a modified SIMP model introduces a small minimum Young's modulus $E_{\min}$:
\begin{equation*}
	E(\rho) = E_{\min} + \rho^p(E_0 - E_{\min}),\quad \rho \in [0,1],
\end{equation*}
where $E_{\min} = 10^{-9} E_0$. This modification ensures that the stiffness matrix remains positive definite at $\rho = 0$, preventing numerical failure \cite{sigmund2007morphology}.

However, optimization results may still exhibit intermediate densities, or "gray areas", which can complicate manufacturing. To address this, techniques like Heaviside projection and sensitivity/density filtering are commonly used to achieve clearer topologies and mitigate numerical issues such as checkerboard patterns.

\subsection{Compliance Minimization Problem}
The goal of compliance minimization with a volume constraint is to find a material density field $\rho(x)\in[0,1]$ that minimizes the total strain energy of an elastic body occupying a bounded domain $\Omega\subset\mathbb{R}^d~(d=2,3)$ under prescribed loads and boundary conditions, while the available material does not exceed a given volume $V^*$. Throughout this paper we assume the displacement $u\in{H}^1(\Omega)$, the density $\rho\in{L}^2(\Omega)$ and the material stiffness tensor $D(\rho)\in{L}^{\infty}(\Omega)$.

The compliance minimization problem in its continuous form can be formulated as follows:
\begin{equation}\label{eq:cmpc}
	\begin{aligned}
		\min_{\rho}:
		&~c(\rho) = \int_{\Omega}\sigma(u):\varepsilon(u)~\mathrm{d}x \\
		\text{subject~to}:
		&~g(\rho) = v(\rho) - V^*\leq0,\quad0\leq\rho(x)\leq1,\\
		&~-\nabla\cdot\sigma(u) = f\quad\text{in}~\Omega,\\
		&~u = 0~\text{on}~\Gamma_D,\quad\sigma(u)\cdot{n} = t~\text{on}~\Gamma_N,
	\end{aligned}
\end{equation}
where $c(\rho)$ is defined as $\int_{\Omega}\sigma(u):\varepsilon(u)~\mathrm{d}x$, this definition equals twice the strain energy, a convention widely adopted in TO literature, $g(\rho) = v(\rho) - V^*\leq0$ denotes the volume constraint, $v(\rho) = \int_{\Omega}\rho~\mathrm{d}x$ is the material volume, $\varepsilon(u) = \frac{1}{2}(\nabla{u}+\nabla{u}^T)$ is the strain tensor, and $\sigma(u) = D(\rho):\varepsilon(u)$ is the stress tensor. Here $f$ and $t$ denote body forces and surface tractions, respectively, and $\Gamma_D \cup \Gamma_N = \partial\Omega$.

Discretising $\Omega$ into $N_e$ finite elements and letting $\bs{\rho} = [\rho_1, \rho_2, \dots, \rho_{N_e}]^T$ be the element-wise constant densities, the problem becomes
\begin{equation}\label{eq:cmpd}
	\begin{aligned} 
		\min_{\bs{\rho}}:
		&~c(\bs{\rho}) = \mathbf{F}^T\mathbf{U}(\bs{\rho})\\ 
		\text{subject~to}:
		&~g(\bs{\rho}) = v(\bs{\rho})- V^*\leq0,\quad0\leq\rho_e\leq1,\\
		&~\mathbf{K}(\bs{\rho})\mathbf{U} = \mathbf{F},\\ 
	\end{aligned}
\end{equation}
with
\begin{equation*}
	v(\bs{\rho}) = \sum_{e=1}^{N_e} v_e\rho_e,\quad\mathbf{K}(\bs{\rho}) = \sum_{e=1}^{N_e}E(\rho_e)\mathbf{K}_e^0.
\end{equation*}
Here $v_e$ is the volume of element $e$, $E(\rho_e)$ follows the SIMP interpolation introduced in Section~\ref{sec:DbM}, and $\mathbf{K}_e^0$ is the elemental matrix for unit Young's modulus.

Solving $\mathbf{K}(\bs{\rho})\mathbf{U}(\bs{\rho}) = \mathbf{F}$ yields the discrete displacement vector $\mathbf{U}$, after which the objective $c(\bs{\rho})$ and its sensitivities are evaluated for the optimization algorithm.

\subsection{Optimization Algorithms}
The Optimality Criteria (OC) method is a classical TO algorithm widely used for compliance minimization problems with volume constraints. It iteratively adjusts the material density $\rho_e$ to satisfy the optimality conditions, based on the update factor:
\begin{equation*}
	B_e = -\frac{\partial{c}(\rho)}{\partial\rho_e}\left(\lambda\frac{\partial{g}(\rho)}{\partial\rho_e}\right)^{-1},
\end{equation*}
where $\lambda$ is the Lagrange multiplier associated with the volume constraint.
 
Bendsøe~\cite{bendsoe1995optimization} proposed a heuristic update scheme for the design variables:
\begin{equation*}
	\rho_e^{\text{new}}=
	\begin{cases}
		\max(0,~\rho_e-m)\quad&\text{if}\quad \rho_eB_e^\eta\leq\max(0,~\rho_e-m)\\
		\min(1,~\rho_e+m)\quad&\text{if}\quad\rho_eB_e^\eta \geq \min(1,~\rho_e+m)\\
		\rho_eB_e^\eta\quad&\text{if}\quad\text{otherwise}
	\end{cases}
\end{equation*}
where $m$ is the move limit and $\eta\in(0,1]$ is a damping exponent. The algorithm is summarized in Algorithm~\ref{alg:oc}.
\begin{algorithm}[htbp]
	\caption{OC pseudo-code}
	\label{alg:oc}
	\begin{minipage}{\textwidth}
		\KwIn{Initial design variables $\bs{\rho}^{(0)}$, volume constraint $V^*$, move limit $m$, damping factor $\eta$, maximum iterations $\text{MaxIter}$, tolerance $\epsilon$}
		\KwOut{Optimized design variables \(\bs{\rho}\)}
		Set \(\bs{\rho}^{(k)} \leftarrow \bs{\rho}^{(0)}\), \(k \leftarrow 0\);
		
		\While{\(k < \mathrm{MaxIter}\) and \(\|\bs{\rho}^{(k+1)} - \bs{\rho}^{(k)}\|_\infty > \epsilon\)}
		{
			Compute objective, constraint, and sensitivities;
			
			(Optional) Apply filters;
			
			Update \(\lambda\) using bisection;
			
			Update \(\bs{\rho}^{\text{new}}\) using optimality criterion;
			
			Set \(\bs{\rho}^{(k+1)} \leftarrow \bs{\rho}^{\text{new}}\), \(k \leftarrow k + 1\);
		}
	\end{minipage}
\end{algorithm}

In TO, in addition to the OC method, the Method of Moving Asymptotes (MMA) is a widely used gradient-based optimization algorithm proposed by Svanberg \cite{svanberg1987method}. It is particularly well-suited for problems involving multiple constraints or complex objective functions. The core idea of MMA is to dynamically adjust the approximation range of the design variables using "moving asymptotes", thereby transforming the original nonlinear optimization problem into a sequence of convex subproblems. This strategy ensures numerical stability throughout the iterative process.

In the compliance minimization problem, MMA approximates the original problem by the following convex subproblem:
\begin{equation*}
	\begin{aligned}
		\min:
		&~\tilde{f}_0^{(k)}(\boldsymbol{\rho}) + a_0z + \sum_{i=1}^{m_c}(c_iy_i+\frac{1}{2}d_iy_i^2)\\
		\text{subject~to}:
		&~\tilde{f}_i^{(k)}(\boldsymbol{\rho}) - a_iz - y_i \leq 0,\quad&i=1,\cdots,m_c\\
		&~\alpha_j^{(k)}\leq{\rho}_j\leq\beta_j^{(k)},\quad&j=1,\cdots,n\\
		&~{y}_i\geq0,~z\geq0,
	\end{aligned}
\end{equation*}
where $m_c$ is the number of constraints, $n$ is the number of design variables, $y_i$ are auxiliary variables for the constraints, and $z$ is an auxiliary variable for the objective. The approximation $\tilde{f}_0^{(k)}(\bs{\rho})$ is the convex approximation of the objective function, while $\tilde{f}_i^{(k)}(\bs{\rho})$ approximates the i-th constraint, both dynamically adjusted based on current gradient information. The bounds $\alpha_j^{(k)}$ and $\beta_j^{(k)}$ are the dynamic lower and upper limits for the design variable $\rho_j$ in the k-th iteration. The parameters $a_0, a_i, c_i, d_i$ are given constants. 

The procedure of the MMA method is summarized in Algorithm~\ref{alg:mma}.
\begin{algorithm}[htbp]
	\caption{MMA pseudo-code}
	\label{alg:mma}
	\begin{minipage}{\textwidth}
		\KwIn{Initial design variables \(\bs{\rho}^{(0)}\), problem parameters, MMA parameters}
		\KwOut{Optimized design variables \(\bs{\rho}\)}
		Set \(\bs{\rho}^{(k)} \leftarrow \bs{\rho}^{(0)}\), \(k \leftarrow 0\);
		
		\While{\(k < \mathrm{MaxIter}\) and \(\|\bs{\rho}^{(k+1)} - \bs{\rho}^{(k)}\|_\infty > \epsilon\)}{
			Compute sensitivities;
			
			(Optional) Apply filters;
			
			Update asymptotes and variable bounds;
			
			Construct and solve the MMA subproblem to obtain \(\bs{\rho}^{\text{new}}\);
			
			Set \(\bs{\rho}^{(k+1)} \leftarrow \bs{\rho}^{\text{new}}\), \(k \leftarrow k + 1\)\;
		}
	\end{minipage}
\end{algorithm}

\subsection{Filtering Methods}
In TO, numerical issues such as mesh dependency, checkerboard patterns, and local minima often arise~\cite{bendsoe2013topology}. Filtering techniques, as a form of regularization, smooth the design variables or sensitivities to mitigate these problems and enhance the physical realizability of optimized structures~\cite{sigmund1998numerical,sigmund2007morphology}.

The sensitivity filter, proposed by Sigmund~\cite{sigmund1997design}, smooths the sensitivities via a weighted average:
\begin{equation*}
	\widetilde{\frac{\partial{f}}{\partial\rho_i}} = \frac{1}{\max\{\gamma,\rho_i\}\sum_{j\in{N}_i}H_{ij}}
	\sum_{j\in{N}_i}H_{ij}\rho_j\frac{\partial{f}}{\partial\rho_j},
\end{equation*}
where $H_{ij} = \max\{0, r_{\min}-\text{dist}(i,j)\}$ is a linearly decaying weight, $N_i$ is the set of elements within filter radius $e_{\min}$, and $\gamma = 10^{-3}$ is a stabilizing constant.

The density filter, proposed by Bruns and Tortorelli~\cite{bruns2001topology} and mathematically justified by Bourdin~\cite{bourdin2001filters}, smooths the material distribution:
\begin{equation*}
	\tilde\rho_i = \frac{\sum_{j\in{N}_i}H_{ij}v_j\rho_j}{\sum_{j\in{N}_i}H_{ij}v_j},
\end{equation*}
where $\tilde{\rho}_i$ is the filtered physical density, $\rho_j$ is the original design variable associated with element $j$, and $v_j$ is the volume of element $j$.

The density filter computes a weighted average of densities using weights $H_{ij}$, smoothing the material distribution. The filtered density $\tilde{\rho}_i$, rather than the design variable $\rho_i$, is used to evaluate structural performance and constraints, serving as the final design in practice.

Under density filtering, the sensitivity of a function $\psi$ with respect ot $\rho_j$ is:
\begin{equation*}
	\frac{\partial\psi}{\partial{\rho}_j} = \sum_{i\in{N}_j}\frac{\partial\psi}{\partial\tilde\rho_i}\frac{\partial\tilde\rho_i}{\partial\rho_j} = \sum_{i\in{N}_j}\frac{H_{ij}v_j}{\sum_{k\in{N}_i}H_{ik}v_k}\frac{\partial\psi}{\partial\tilde\rho_i},
\end{equation*}
where $\psi$ is an objective or constraint function, and $N_j$ s the set of elements influenced by $\rho_j$.

Heaviside projection filtering~\cite{guest2004achieving} projects the filtered density to a physical density $\bar{\rho}_i$, achieving a clear black-and-white topology and enforcing a minimum length scale. The projection is defined as:
\begin{equation*}
	\bar\rho_i= 1 - e^{-\beta\tilde\rho_i} + \tilde\rho_ie^{-\beta},
\end{equation*}
where $\beta$ controls smoothness. A continuation scheme gradually increases $\beta$ to enhance the projection effect and ensure numerical stability.

After applying the Heaviside projection filter, the sensitivities of a function with respect to $\tilde\rho_i$ must also be computed using the chain rule:
\begin{equation*}
	\frac{\partial\psi}{\partial\tilde\rho_i} = \frac{\partial\psi}{\partial\bar\rho_i}(\beta{e}^{-\beta\tilde\rho_i} + e^{-\beta}).
\end{equation*}

The filter radius $r_{\min}$ significantly affects the result: small $r_{\min}$ leads to oscillations or checkerboard patterns, while large $r_{\min}$ causes excessive smoothing and loss of design details.

\section{SOPTX Framework Design and Multi-Backend Switching}
This section introduces the design principles and key components of the SOPTX framework, emphasizing its multi-backend switching capability for improved computational performance and flexibility on central processing units (CPUs) and graphics processing units (GPUs). The section is divided into three parts: an overview of the FEALPy architecture as the tensor computation foundation, a description of SOPTX's modular structure (material, solver, filter, and optimization modules), and an analysis of the multi-backend mechanism supporting NumPy, PyTorch, and JAX for diverse computational needs.

\subsection{FEALPy Architecture Design}
\begin{figure}[!htbp]
	\centering
	\includegraphics[width=0.85\textwidth]{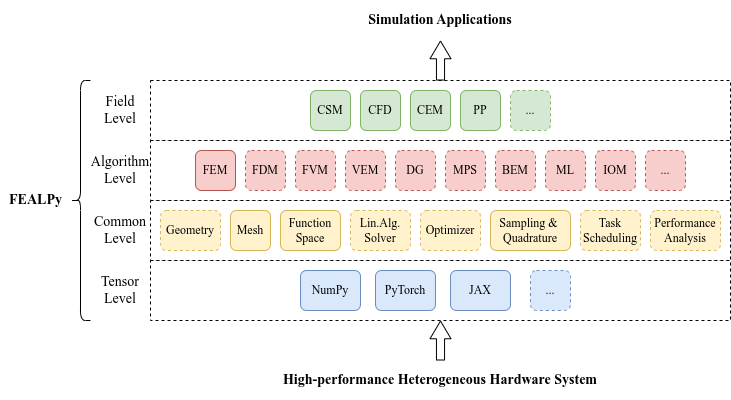}
	\caption{The layered architecture of FEALPy, comprising tensor, common, algorithm, and field levels, progressing from low-level functionalities to high-level applications. Modules in dashed boxes are under development.}
	\label{fs:fig1}
\end{figure}

As shown in Figure~\ref{fs:fig1}, FEALPy adopts a layered architecture with four levels, arranged from bottom to top: tensor, common, algorithm, and field. Building on this, FEALPy introduces the \textit{Tensor Backend Manager}, which unifies the management of computational backends like NumPy, PyTorch, and JAX. This manager provides a consistent tensor operation interface, following the Python Array API Standard v2023.12~\cite{arrayapi2023}, enabling seamless adaptation across various software and hardware platforms.

Specifically, the functions of each layer are as follows:
\begin{enumerate} 
	\item \textbf{Tensor level}: Provides core tensor operations and manages backends via the \textit{Tensor Backend Manager}, supporting SOPTX's multi-backend switching. It also enables automatic differentiation (AD) in PyTorch and JAX, simplifying sensitivity computations in topology optimization (TO).
	\item \textbf{Common level}: Includes mesh generation and finite element spaces, supporting rapid construction of meshes and function spaces for finite element analysis in SOPTX.
	\item \textbf{Algorithm level}: Encompasses solvers and optimization algorithms, offering efficient and stable computational support for SOPTX's optimization methods.
	\item \textbf{Field level}: Targets specific physical problems (e.g., linear elasticity), enabling SOPTX to address diverse TO problems, such as compliance minimization under volume constraints, with tailored application support.
\end{enumerate}

This layered design ensures FEALPy's functionality, extensibility, and performance, with the \textit{Tensor Backend Manager} abstracting backend differences to enhance user focus on algorithms and applications.

\subsection{SOPTX Architecture Design}
The SOPTX framework is positioned within the field level of FEALPy's layered architecture, targeting structural topology optimization (STO) applications. SOPTX fully inherits and extends FEALPy's \textit{Tensor Backend Manager}. It also leverages a variety of numerical algorithm components and generic mesh and geometry handling capabilities. This design ensures flexibility and extensibility while enhancing computational efficiency for TO problems.

\begin{figure}[!htb]
	\centering
	\includegraphics[width=0.85\textwidth]{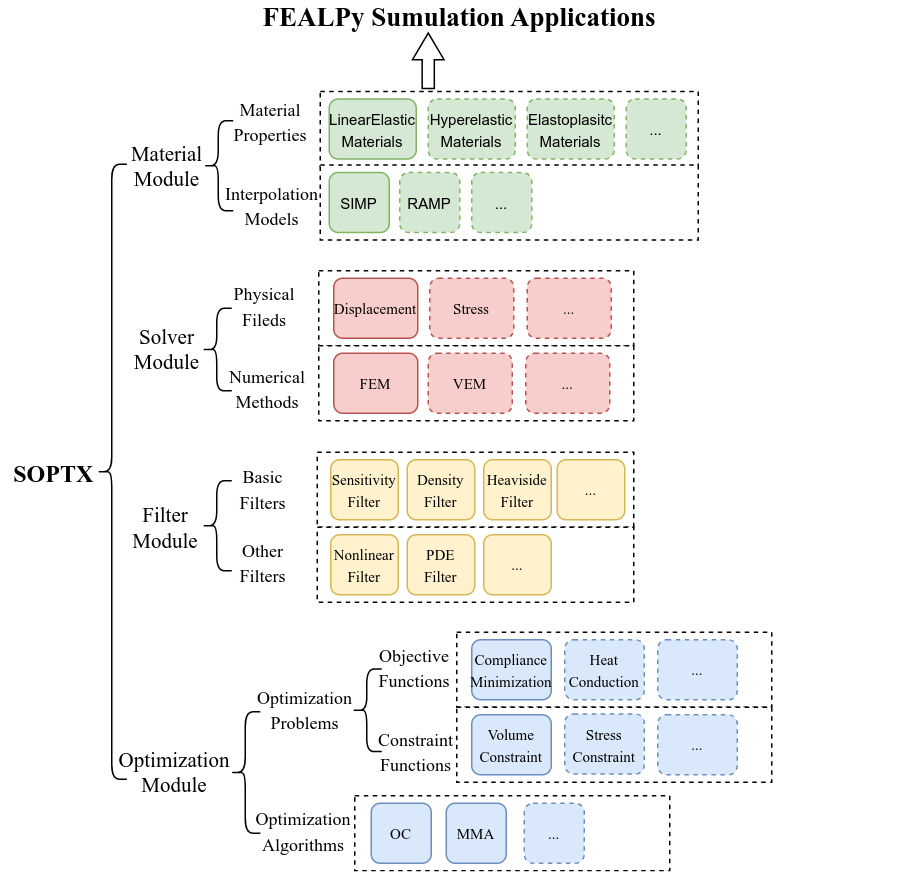}
	\caption{The modular architecture of SOPTX, consisting of material, solver, filter, and optimization modules. The material module serves as the foundation, while the solver and filter modules manage intermediate computations, collectively supporting the optimization module. Components in dashed boxes are under development.}
	\label{fs:fig2}
\end{figure}

As shown in Figure~\ref{fs:fig2}, SOPTX adopts a modular architecture consisting of four primary components: the material, solver, filter, and optimizer modules. These components communicate and share data through clearly defined interfaces, forming a loosely coupled and easily extensible multi-backend framework for TO.

\subsubsection{Material Module}
The material module in SOPTX currently focuses on linear elastic materials, extending FEALPy's implementation with TO-specific interfaces for computing Lamé constants, elasticity matrices, and strain matrices. It implements the Solid Isotropic Material with Penalization (SIMP) model, commonly employing a penalization factor of $p=3$, to penalize  continuous density fields towards discrete material distributions. The module's abstract interface design ensures extensibility. This allows future support for other models, such as the Rational Approximation of Material Properties (RAMP), or complex behaviors like anisotropy, hyperelasticity, and elastoplasticity through subclassing.

The module encapsulates the interpolation, update, and evaluation of elastic constants within a unified interface, providing essential material data for downstream modules. Its key functional interfaces include:
\begin{enumerate}
	\item Computing elasticity tensors from density fields for the solver module.
	\item Efficient batch updates of material properties for iterative optimization.
	\item Providing base material parameters (e.g., elastic constants) for optimization.
\end{enumerate}

These interfaces ensure consistent backend compatibility and optimized performance. The material module bridges physical modeling and numerical optimization, maintaining independence while supplying essential data throughout the TO workflow. Its modular design simplifies current linear elasticity solutions and provides a flexible foundation for future material model extensions.

\subsubsection{Solver Module}
The solver module forms the computational core of the SOPTX framework, tasked with solving physical field problems (e.g., displacement fields) based on given material properties, boundary conditions, and external loads. The module is designed to be efficient, flexible, and backend-independent. This design addresses the demands of repeatedly solving large-scale linear systems in TO. To achieve this, the module enhances computational efficiency through optimized matrix assembly, intelligent caching, and multi-backend acceleration, supporting backends including NumPy, PyTorch, and JAX.

In its current version, SOPTX primarily adopts the finite element method (FEM) to solve linear elasticity equations, leveraging FEALPy's finite element module. The key features of the solver module include:
\begin{itemize}
	\item Dimensional and element adaptability: It supports various spatial dimensions, element types, and boundary conditions, where the density field is represented as piecewise constant per element and the displacement field uses linear finite elements, ensuring both stability and computational efficiency.
	\item Efficient numerical strategies: The module implements a fast matrix assembly technique that separates element-independent and element-dependent components to eliminate redundant computations. Additionally, it incorporates symbolic integration to boost efficiency by precomputing exact expressions.
	\item Multiple solution strategies: It provides both direct solvers (e.g., MUMPS) and iterative solvers (e.g., Conjugate Gradient, CG), with automatic optimization tailored to the backend, such as utilizing GPU acceleration on PyTorch and JAX.
\end{itemize}

The solver module is designed with strong extensibility, enabling the future incorporation of additional physical fields (e.g., stress fields) and numerical methods (e.g., Virtual Element Method), as well as adaptation to nonlinear mechanics and multiphysics coupling scenarios."

Within the SOPTX framework, the solver module acts as a bridge between physical analysis and optimization computation through the following interactions:
\begin{itemize}
	\item Interaction with the material module: It receives material properties, such as elasticity matrices, via standardized interfaces, ensuring the solver remains independent of specific material interpolation schemes.
	\item Output to the optimization module: It provides the displacement field for evaluating the objective function and the stiffness matrix for sensitivity analysis.
	\item Result reuse: It employs an intelligent caching mechanism to avoid recomputing invariant components. This approach significantly reduces computational overhead.
\end{itemize}

Through its modular design and well-defined interfaces, the solver module efficiently conducts physical field computations, providing critical support for the TO process while ensuring both consistency and high performance across multiple backends.

\subsubsection{Filter Module}
The filter module in SOPTX is a cornerstone of TO, tasked with processing design variables and applying regularization techniques. Its primary functions are to mitigate numerical instabilities, such as checkerboarding and mesh dependency, and to improve the manufacturability of optimized structures. Adhering to the framework's loosely coupled design, the filter module operates as an independent unit with efficient data exchange across modules. It offers a unified interface that supports multiple computational backends. Backend-specific optimizations ensure high computational efficiency.

The module implements three key filtering techniques:
\begin{itemize}
	\item Sensitivity filtering: Applies weighted averaging to objective function sensitivities, which effectively suppresses checkerboard patterns.
	\item Density filtering: Maps raw design variables to physical densities, addressing the locality limitations of sensitivity filtering.
	\item Heaviside projection filtering: Enhances density filtering with a smoothed Heaviside function to achieve distinct black-and-white designs by driving intermediate densities toward 0 or 1.
\end{itemize}

The module employs a KD-tree-based neighborhood search to construct filtering matrices, which enables efficient regularization on unstructured meshes and complex geometries. Optimization strategies, such as precomputed neighborhoods and sparse matrix representations, ensure scalability for large-scale TO problems.

The filter module is highly extensible. Users can subclass the base filter class to implement custom algorithms, such as nonlinear filters or PDE-based filters. This approach seamlessly inherits multi-backend compatibility. Additionally, the framework supports chained filter combinations, allowing for customized regularization strategies to meet various optimization needs.

The filter module interacts with other components in the following ways:
\begin{itemize}
	\item Processing design variables: It processes raw design variables into filtered physical densities.
	\item Interaction with optimization module: It filters unprocessed sensitivities ensuring numerical stability during optimization.
	\item Interaction with material module: It provides filtered densities for material interpolation, defining a streamlined data flow:
	\begin{equation*}
		\text{design variables} \to \text{filtering} \to \text{physical density} \to \text{material properties} \to \text{physical response}
	\end{equation*}
\end{itemize}

Through its modular and extensible design, the filter module not only addresses numerical challenges in TO but also lays a robust foundation for advanced regularization, enabling high-quality, manufacturable structural designs within the SOPTX framework.

\subsubsection{Optimization Module}
The optimization module serves as the computational core of the SOPTX framework, tasked with defining and solving optimization problems by integrating physical field solutions with user-defined objectives. It collaborates seamlessly with the material, solver, and filter modules to drive the TO pipeline. Adhering to the framework's multi-backend design, it offers a unified interface across NumPy, PyTorch, and JAX, utilizing backend-specific optimizations, including GPU acceleration in PyTorch and JAX, to enhance performance.

The module currently focuses on compliance minimization under volume constraints, a foundational problem in TO. It supports two mainstream algorithms: Optimality Criteria (OC) and the Method of Moving Asymptotes (MMA). OC is a lightweight method for volume-constrained problems, while MMA, based on Svanberg's 2007 formulation~\cite{Svanberg2007MmaAG}, has been adapted for multi-backend compatibility and optimized for efficiency.

The design of the optimization module follows two core principles
\begin{enumerate}
	\item Decoupling of problem and algorithm: The optimization problem (e.g., objectives, constraints, sensitivities) is separated from the solving algorithm, which enables flexible pairing of formulations and solvers.
	\item Dual-mode sensitivity: Supports both manually derived sensitivities for transparency and AD in PyTorch or JAX for ease and scalability.
\end{enumerate}

This design ensures high extensibility, enabling users to define custom problems by subclassing the base problem class to introduce new objectives (e.g., heat conduction, compliant mechanisms) or constraints (e.g., stress, frequency). Additionally, users can integrate other optimization algorithms, such as Sequential Quadratic Programming (SQP) or Sequential Linear Programming (SLP), in a modular fashion.

The optimization module serves as a central hub within the SOPTX framework, interacting with other modules as follows:
\begin{itemize}
	\item Interaction with solver module: It utilizes physical outputs, such as displacement fields and stiffness matrices, to compute objectives and sensitivities.
	\item Interaction with Filter Module: It exchanges sensitivities with the filter module for processing and updates, and receives filtered physical densities derived from design variables.
\end{itemize}

With its modular and extensible architecture, the optimization module effectively addresses classical TO problems while providing a versatile foundation for advanced multiphysics optimization. It acts as the decision-making hub of SOPTX, balancing computational performance, flexibility, and user adaptability.

\subsection{Multi-Backend Switching}
In STO, computational performance is a central concern in both research and engineering. To meet diverse computational demands, the SOPTX framework builds upon FEALPy to implement a flexible multi-backend architecture. This allows seamless switching between tensor computation backends such as NumPy, PyTorch, and JAX, enhancing the software's applicability, efficiency, portability, and flexibility across various hardware and software platforms.

Specifically, each backend offers distinct advantages for different scenarios
\begin{itemize} 
	\item \textbf{NumPy}: Suited for small-scale tasks and rapid prototyping, offering stable performance and efficient memory management on CPU platforms.
	\item \textbf{PyTorch and JAX}: Both offer GPU acceleration and AD, making them ideal for large-scale or high-dimensional problems. AD facilitates sensitivity analysis, while JAX provides just-in-time (JIT) compilation and automatic vectorization, further enhancing performance on GPU platforms.
\end{itemize}

\section{Getting Started with SOPTX}
This section provides installation instructions and usage guidance for SOPTX, enabling readers to quickly grasp the framework's operations and apply it to topology optimization (TO) tasks. Built on FEALPy, SOPTX leverages a multi-backend switching mechanism and modular design to deliver an efficient and flexible computational environment. The section is structured into two sections: first, a detailed guide on installing SOPTX and its dependency FEALPy, ensuring proper configuration of the development environment; second, a demonstration of SOPTX's workflow through a classical 2D cantilever beam compliance minimization example.

\subsection{Software Installation}
SOPTX is a TO toolkit built on top of FEALPy, an intelligent CAE simulation engine that provides numerical computing capabilities. Installing FEALPy is required before SOPTX. For flexibility and the latest features, install both from source. Ensure Git and Python are installed. Use a virtual environment to avoid dependency conflicts.

Core installation steps:
\begin{enumerate} 
	\item Clone the FEALPy repository from GitHub:
\begin{lstlisting}[language=bash]
git clone https://github.com/weihuayi/fealpy.git
\end{lstlisting}
	\item Change into the FEALPy directory and install it in editable mode:
\begin{lstlisting}[language=bash]
cd fealpy
pip install -e . 
\end{lstlisting}
	\item Similarly, install SOPTX:
\begin{lstlisting}[language=bash]
git clone https://github.com:weihuayi/soptx.git
cd soptx
pip install -e . 
\end{lstlisting}
\end{enumerate}

For the complete installation guide, please refer to the official documentation at: \url{https://github.com/weihuayi/fealpy}.

\subsection{Example: 2D Cantilever Beam}\label{sec:exp_canti_beam}
\begin{figure}[htbp]
	\centering
	\includegraphics[width=0.5\textwidth]{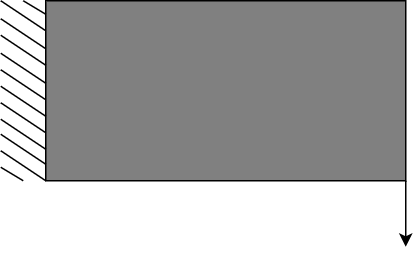}
	\caption{Cantilever beam geometry: fixed on the left, with a downward concentrated load on the right.}
	\label{fig:cantilever}
\end{figure}

We use a popular compliance minimization benchmark to demonstrate the usage of SOPTX: minimizing the structural compliance of a cantilever beam under tip loading (see Figure~\ref{fig:cantilever})~\cite{bendsoe2013topology}. The left end of the beam is fixed, and a downward concentrated load $T = -1$ is applied to the bottom of the right end. A $160 \times 100$ uniform quadrilateral mesh is used. The target volume fraction is set to 0.4, with material properties $E = 1$ and $\nu = 0.3$. The penalization factor is $p = 3$, and the sensitivity filter radius is $r = 6.0$, which matches the mesh element size to ensure structural smoothness and eliminate checkerboard patterns.

To begin, we import essential modules from FEALPy (backend, mesh, function spaces) and SOPTX (material, solver, filter, optimization), as shown in Code~\ref{code:mouldes_pde}. The partial differential equation (PDE) model for the cantilever beam, defined in \texttt{Cantilever2dData1} (see~\ref{sec:code_canti_2d}), specifies the geometry, load, and boundary conditions via standard SOPTX interfaces.

\begin{lstlisting}[language=python, caption={Module imports and PDE model}, label={code:mouldes_pde}]
from fealpy.backend import backend_manager as bm
from fealpy.mesh import UniformMesh
from fealpy.functionspace import LagrangeFESpace, TensorFunctionSpace
from soptx.material import DensityBasedMaterialConfig, DensityBasedMaterialInstance
from soptx.solver import ElasticFEMSolver, AssemblyMethod
from soptx.filter_ import SensitivityBasicFilter
from soptx.opt import ComplianceObjective, ComplianceConfig, VolumeConstraint, VolumeConfig
from soptx.opt import OCOptimizer
from soptx.pde import Cantilever2dData1

pde = Cantilever2dData1(xmin=0, xmax=160, ymin=0, ymax=100,	T = -1)
\end{lstlisting}

Next, we define the mesh and finite element spaces (Code~\ref{code:mesh_space}). The displacement field uses a first-order continuous Lagrange space, while the density field is represented in a zeroth-order discontinuous Lagrange space, aligning with typical TO discretizations.

\begin{lstlisting}[language=python, caption={Mesh and function space definitions}, label={code:mesh_space}]
mesh = UniformMesh(extent=[0, 160, 0, 100], h=[1, 1], origin=[0, 0])
space_C = LagrangeFESpace(mesh=mesh, p=1, ctype='C')
tensor_space_C = TensorFunctionSpace(scalar_space=space_C, shape=(-1, 2))
space_D = LagrangeFESpace(mesh=mesh, p=0, ctype='D')
\end{lstlisting}

The material module is then instantiated with the specified properties and Solid Isotropic Material with Penalization (SIMP) interpolation model (Code~\ref{code:material}). This module handles material property computations based on the density field.

\begin{lstlisting}[language=python, caption={Material module}, label={code:material}]
material_config = DensityBasedMaterialConfig(
		elastic_modulus=1.0, minimal_modulus=1e-9, 
		poisson_ratio=0.3, plane_assumption="plane_stress",  
		interpolation_model="SIMP", penalty_factor=3.0)
materials = DensityBasedMaterialInstance(config=material_config)
\end{lstlisting}

Subsequently, the solver module is initialized with the material module, PDE model, standard matrix assembly, and a direct solver (MUMPS), as shown in Code~\ref{code:solver_filter}. A sensitivity filter is also set up to regularize the optimization.

\begin{lstlisting}[language=python, caption={Solver and filter module}, label={code:solver_filter}] 
solver = ElasticFEMSolver(materials=materials,
			tensor_space=tensor_space_C, pde=pde,
			assembly_method=AssemblyMethod.STANDARD,
			solver_type='direct',
			solver_params={'solver_type': 'mumps'})
sens_filter = SensitivityBasicFilter(mesh=mesh, rmin=6.0)
\end{lstlisting}

The optimization problem is defined by instantiating the compliance objective and volume constraint classes, followed by configuring the Optimality Criteria (OC) optimizer with a maximum of 200 iterations and a convergence tolerance of 0.01 (Code~\ref{code:optimizer}).

\begin{lstlisting}[language=python, caption={Optimization module}, label={code:optimizer}]
objective = ComplianceObjective(solver=solver)
constraint = VolumeConstraint(solver=solver, volume_fraction=0.4)
optimizer = OCOptimizer(objective=objective,
			constraint=constraint, filter=sens_filter,
			options={'max_iterations': 200, 'tolerance': 0.01})
\end{lstlisting}

Finally, the initial density field is uniformly set to the target volume fraction, and the optimization is executed. Results are saved, and convergence history is plotted (Code~\ref{code:main_post}).

\begin{lstlisting}[language=python, caption={Main program and post-processing}, label={code:main_post}]
if __name__ == "__main__":
	@cartesian
	def density_func(x):
		val = config.init_volume_fraction * bm.ones(x.shape[0], **kwargs)
		return val
	rho = space_D.interpolate(u=density_func)
	rho_opt, history = optimizer.optimize(rho=rho[:])
	
	from soptx.opt import save_optimization_history, plot_optimization_history
	save_optimization_history(mesh, history)
	plot_optimization_history(history)
\end{lstlisting}

\begin{figure}[htb]
	\centering
	\includegraphics[width=1.0\textwidth]{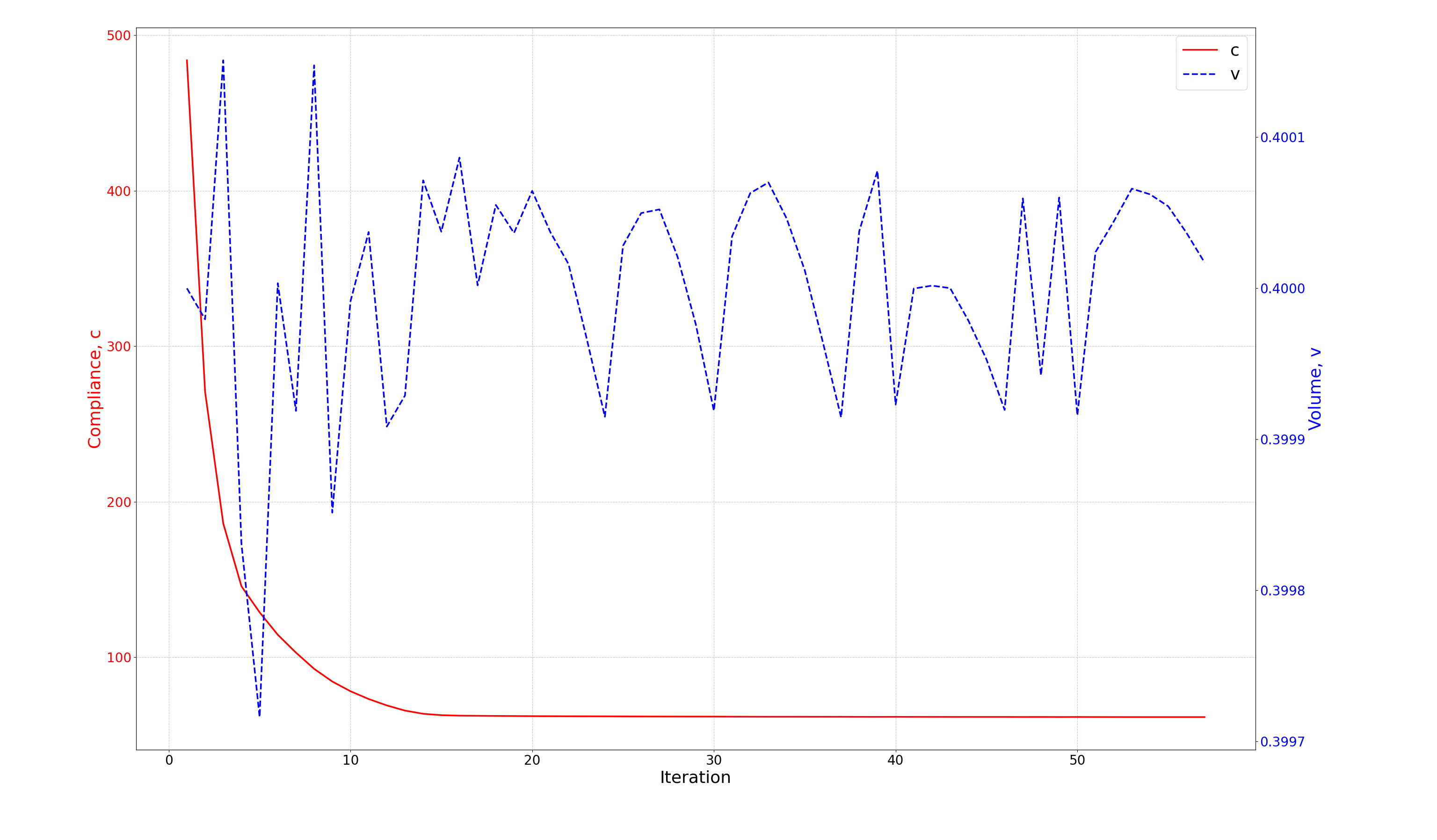}
	\caption{Convergence histories of the compliance $c(\rho)$ and volume fraction $v(\rho)$ for the 2D cantilever beam initialized with a uniform density of 0.4.}
	\label{fig:canti_04_convergence}
\end{figure}

Figure~\ref{fig:canti_04_convergence} shows the convergence histories of the compliance $c(\rho)$ and the volume fraction $v(\rho)$ for an initial density of 0.4. Compliance decreases rapidly from its initial value of approximately 500 to around 100 in the first 10 iterations, converging by iteration 57. The volume fraction remains stably around 0.4, with only minor fluctuations.

\begin{figure}[htb]
	\centering
	\subfloat[$c(\rho)=145.5753,~v(\rho)=0.3998$]{
		\includegraphics[width=0.32\textwidth]{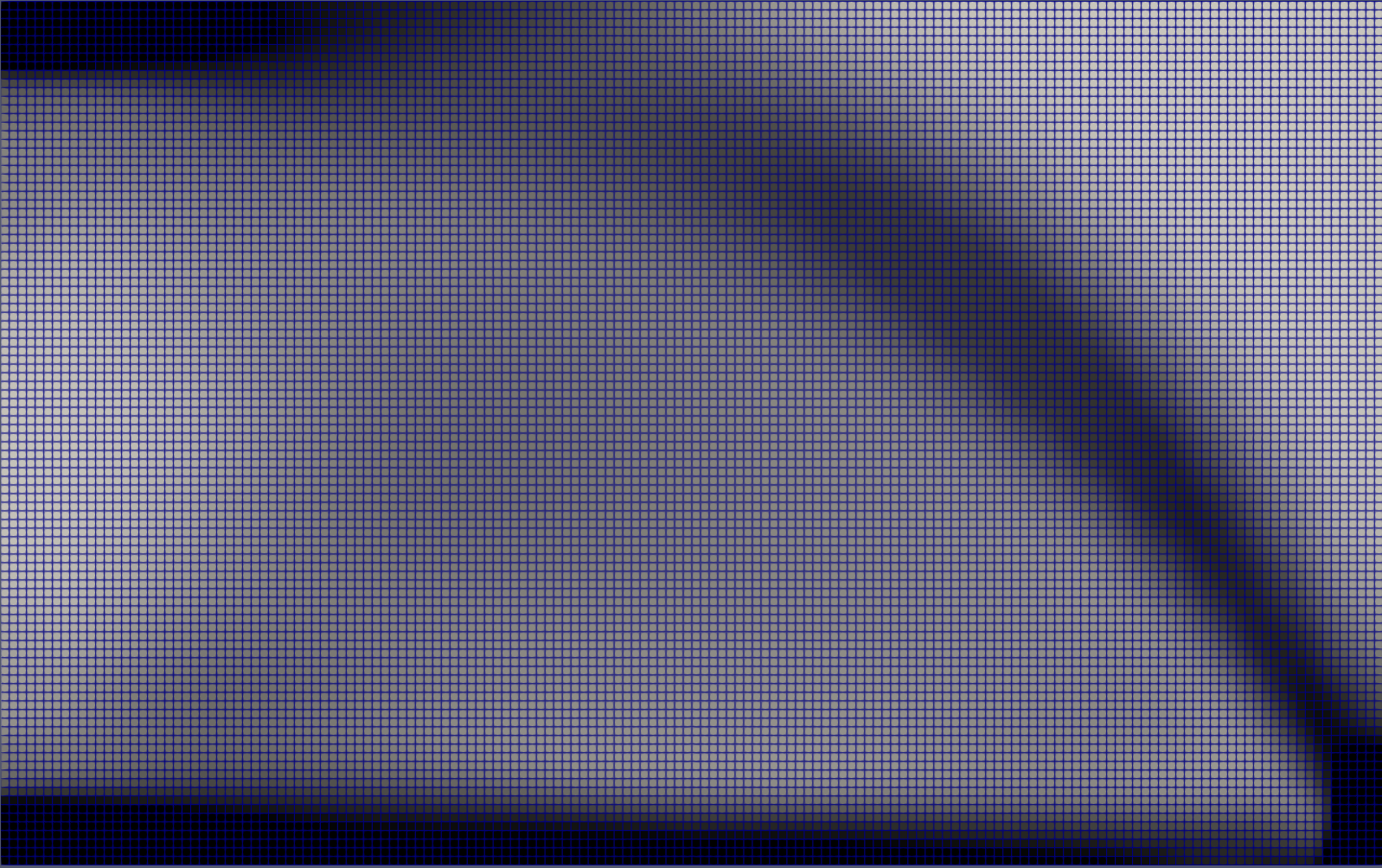}
		\label{fig:canti_04_1}
	}
	\subfloat[$c(\rho)=61.7084,~v(\rho)=0.3999$]{
		\includegraphics[width=0.32\textwidth]{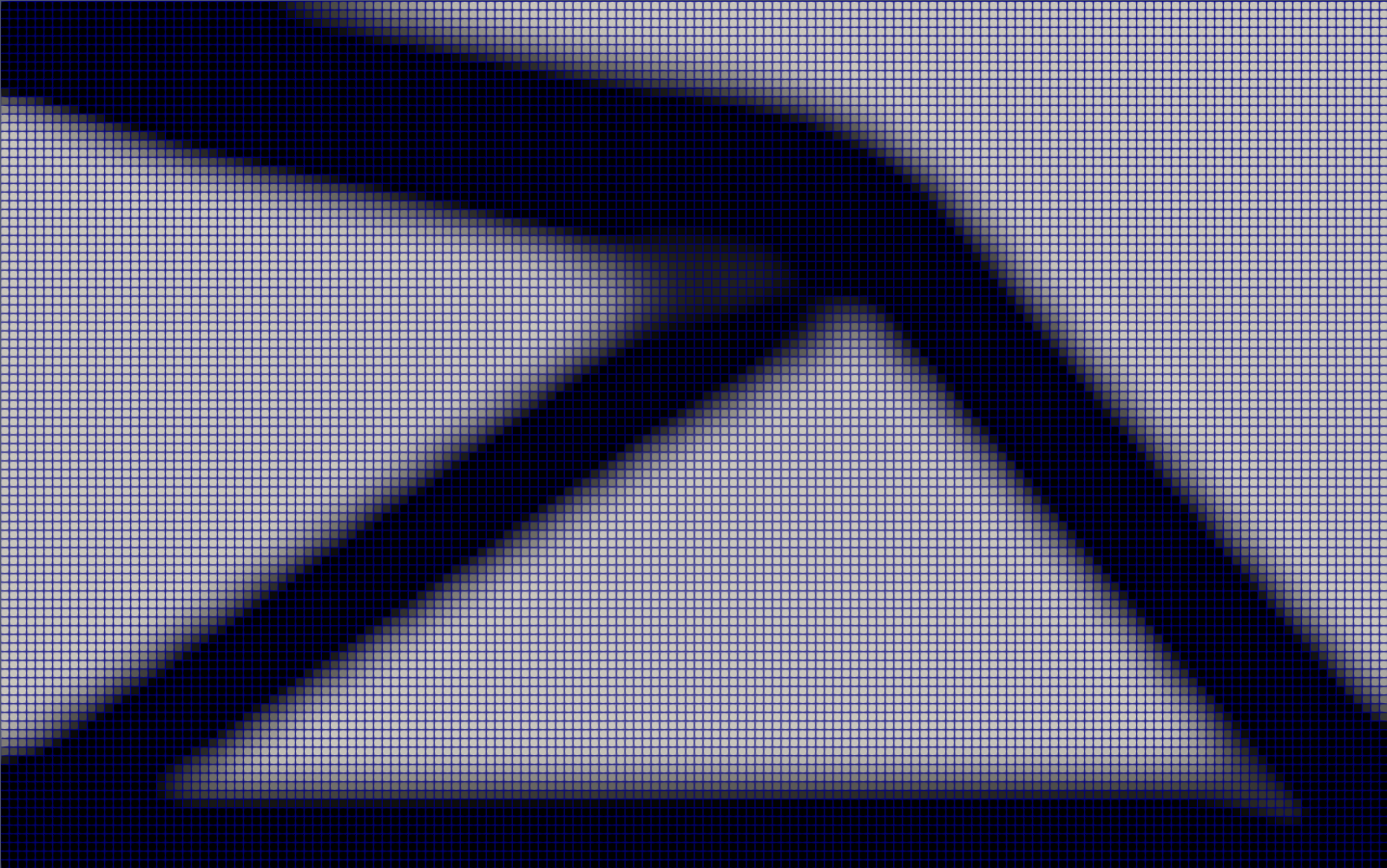}
		\label{fig:canti_04_2}
	}
	\subfloat[$c(\rho)=61.4208,~v(\rho)=0.4000$]{
		\includegraphics[width=0.32\textwidth]{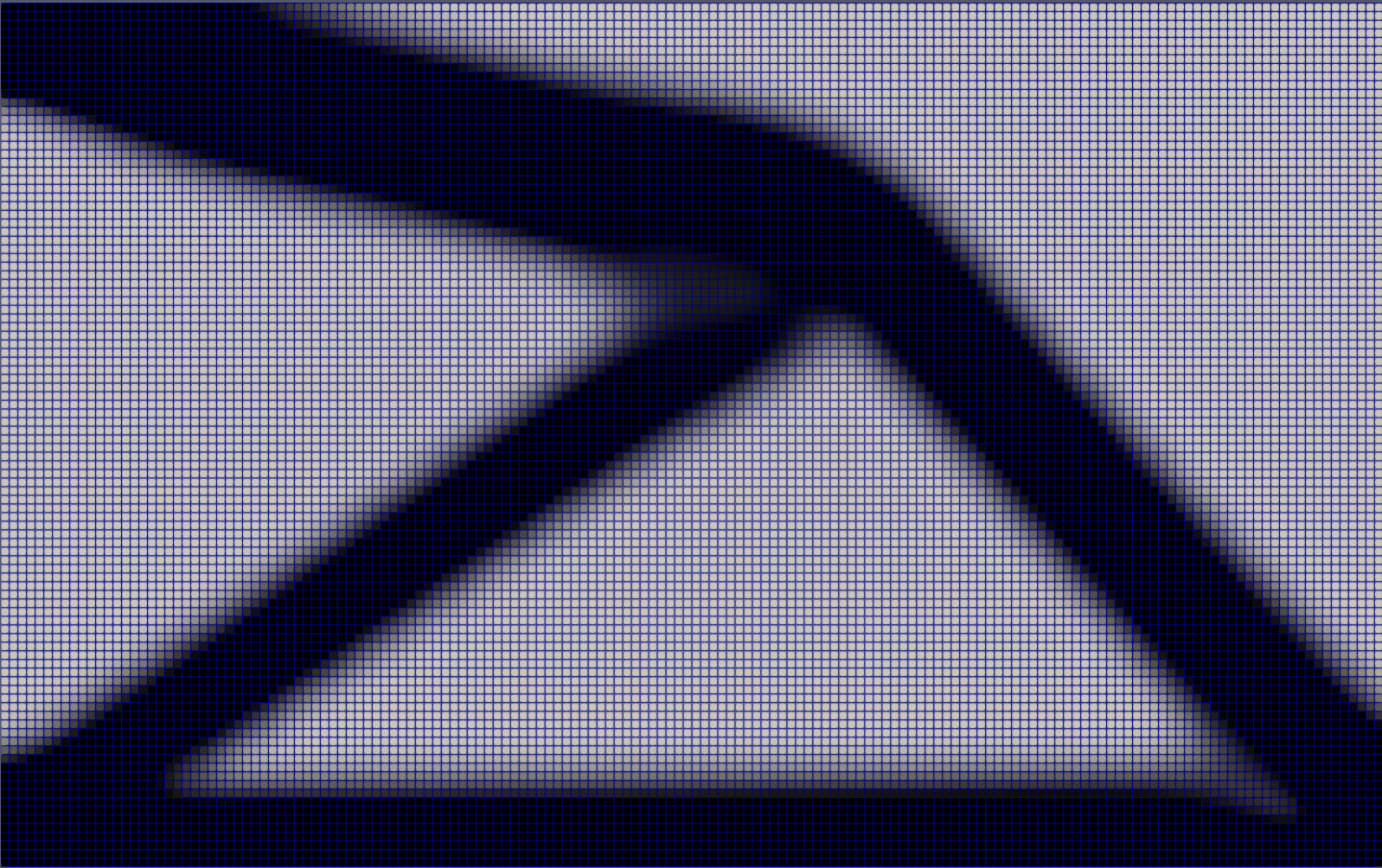}
		\label{fig:canti_04_3}
	}
	\caption{Topology layouts at iterations 3 (left), 30 (middle), and 57 (right) during the optimization process. Each subfigure reports the corresponding compliance and volume fraction values.}
	\label{fig:canti_04_all}
\end{figure}

Figure~\ref{fig:canti_04_all} displays the resulting topologies at iterations 3, 30, and 57.

\begin{figure}[htb]
	\centering
	\includegraphics[width=1.0\textwidth]{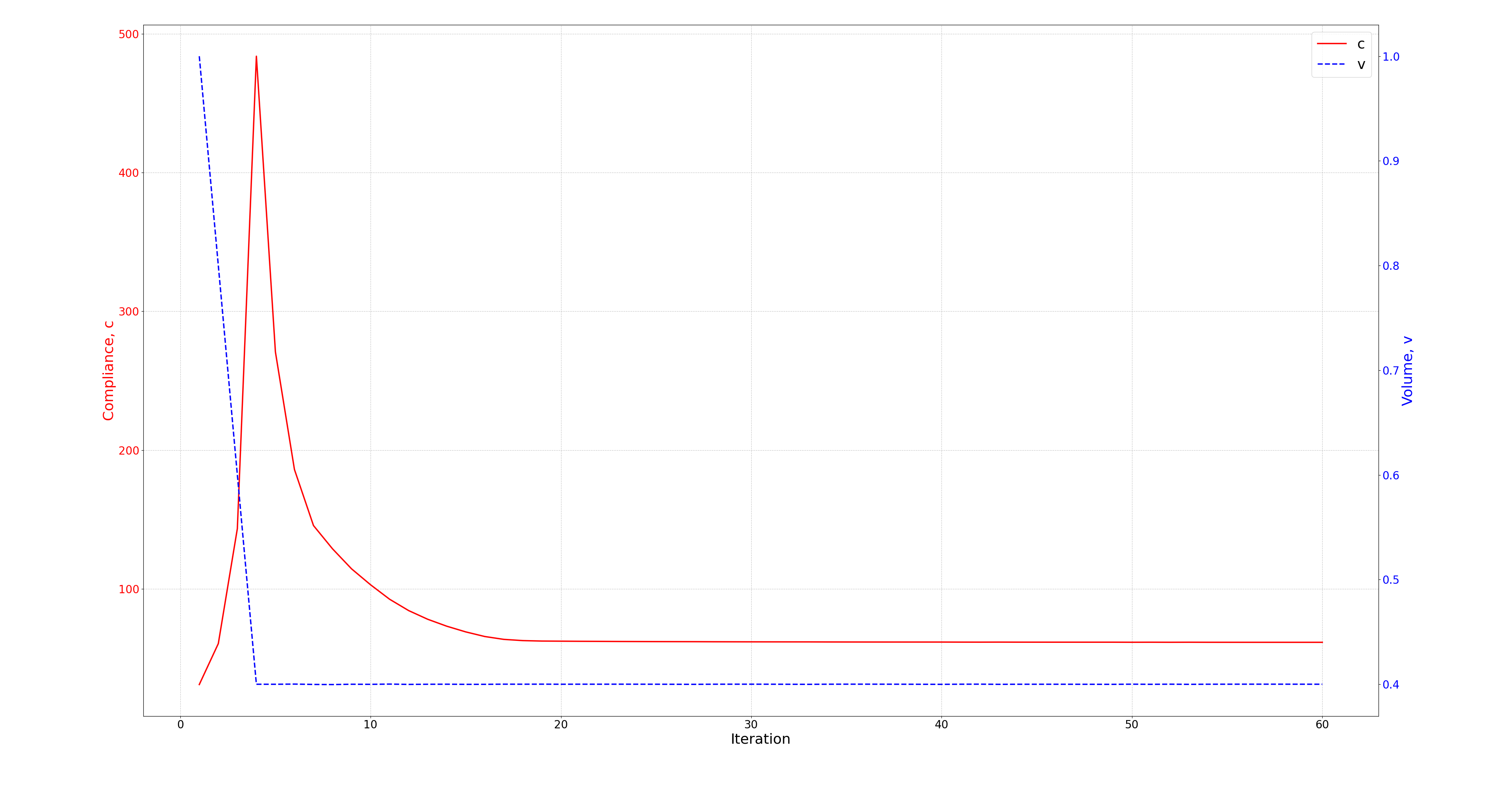}
	\caption{Convergence histories of the compliance $c(\rho)$ and volume fraction $v(\rho)$ for the 2D cantilever beam with an initial density of $1$.}
	\label{fig:canti_1_convergence}
\end{figure}

Figure~\ref{fig:canti_1_convergence} illustrates the convergence behavior for an initial density of 1. The volume fraction decreases from 1 to 0.4 within 5 iterations, while compliance initially increases to around 500 before decreasing to convergence at iteration 60. This process requires slightly more iterations due to the initial adjustment to meet the volume constraint.

\begin{figure}[htb]
	\centering
	\subfloat[$c(\rho)=185.9530,~v(\rho)=0.4002$]{
		\includegraphics[width=0.32\textwidth]{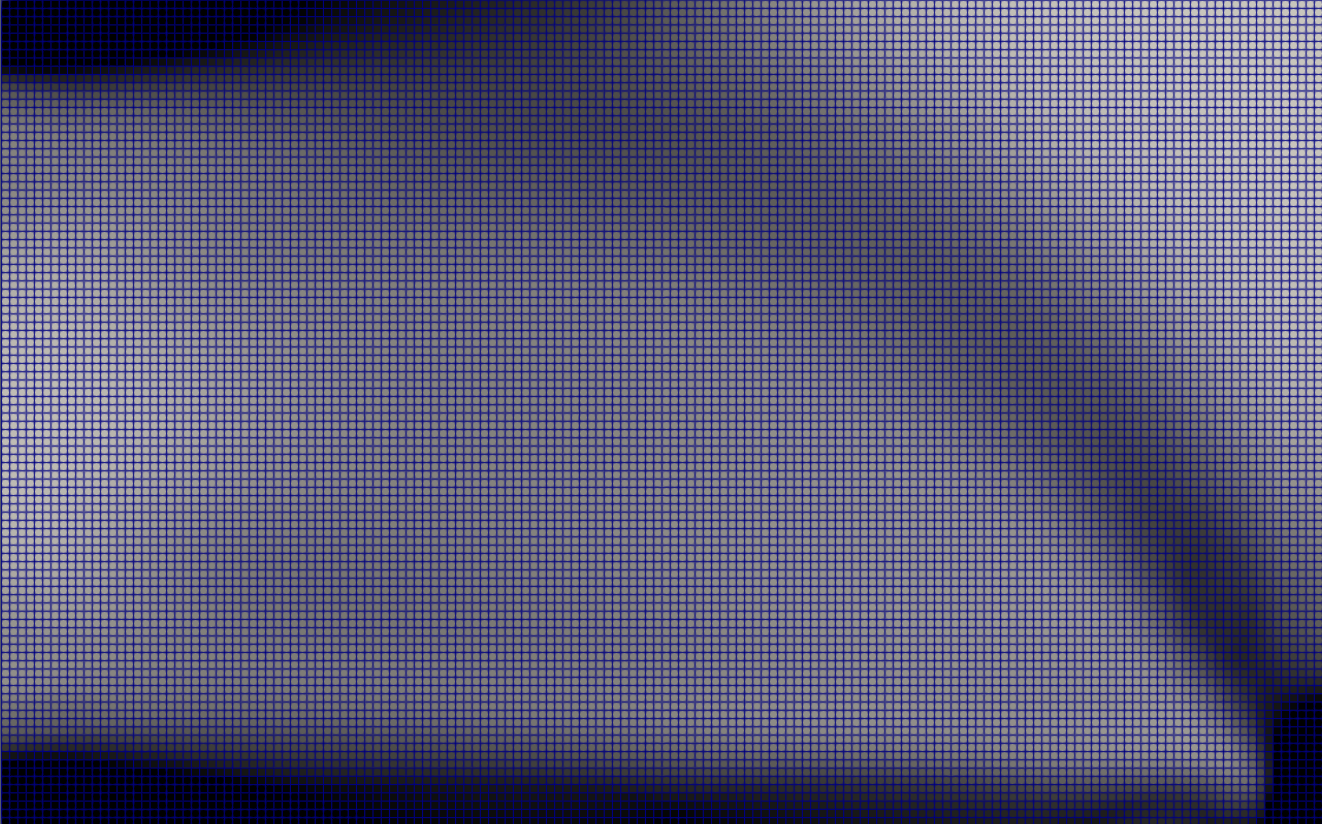}
		\label{fig:canti_1_1}
	}
	\subfloat[$c(\rho)=61.7084,~v(\rho)=0.3999$]{
		\includegraphics[width=0.32\textwidth]{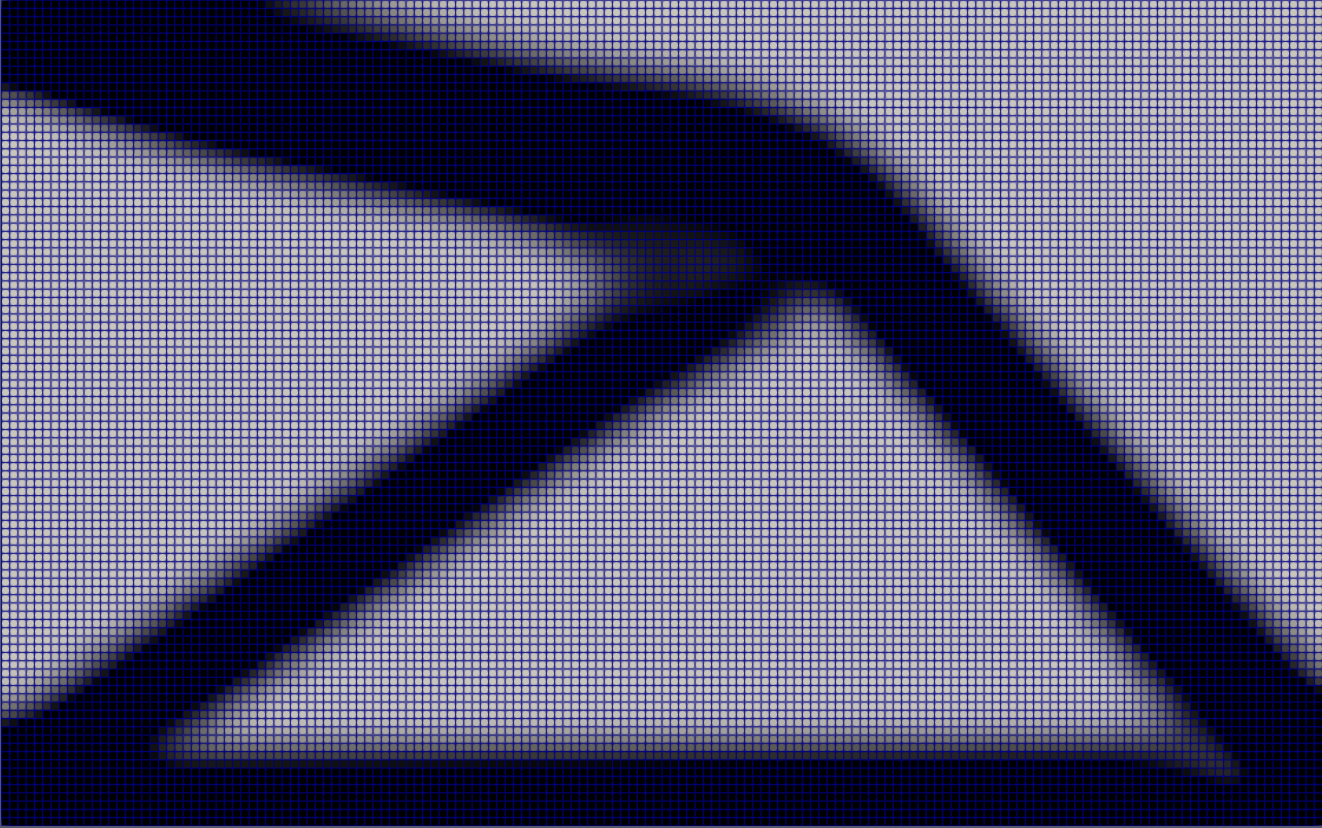}
		\label{fig:canti_1_2}
	}
	\subfloat[$c(\rho)=61.4208,~v(\rho)=0.4000$]{
		\includegraphics[width=0.32\textwidth]{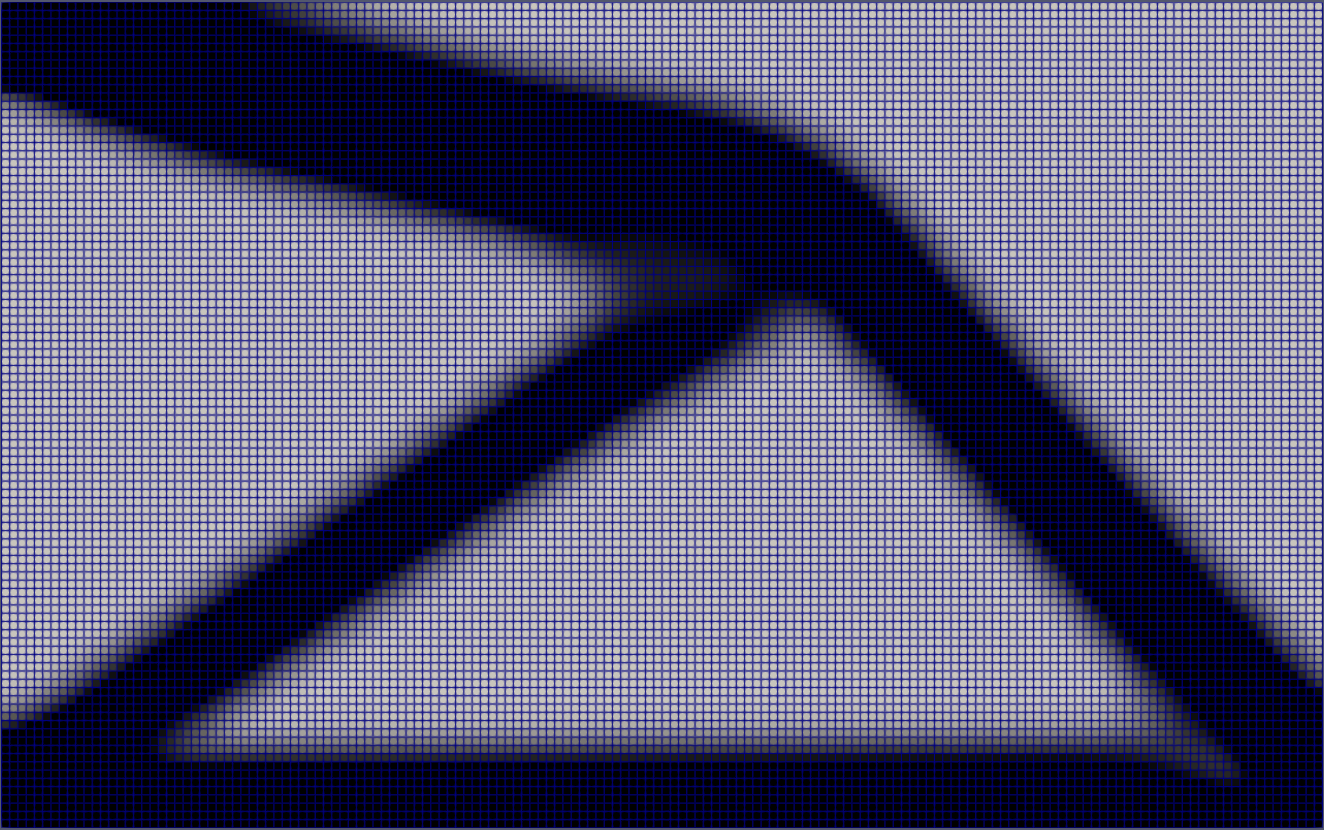}
		\label{fig:canti_1_3}
	}
	\caption{Topology layouts at iterations 6 (left), 33 (middle), and 60 (right) during the optimization process. Each subfigure includes the compliance and volume fraction values.}
	\label{fig:canti_1_all}
\end{figure}

Figure~\ref{fig:canti_1_all} shows the topology layouts at iterations 6, 33, and 60.

These results highlight the robustness of SOPTX in handling different initial conditions and its efficiency in achieving convergence with the OC algorithm.

\section{Numerical Examples}
This section presents numerical examples to validate the SOPTX framework's capabilities in topology optimization (TO), spanning 2D MBB beam problems to 3D cantilever structures. These examples demonstrate SOPTX's versatility across diverse partial differential equation (PDE) models, meshes, filters, and algorithms, while showcasing advanced features like fast matrix assembly, automatic differentiation (AD), and multi-backend switching. All runs are performed on a desktop computer running Ubuntu 24.04, using a CPU AMD Ryzen 9 9950X @ 4.5GHz, 2 × 32 GB of RAM, and an NVIDIA GeForce RTX 5070Ti GPU.

The section is organized as follows: Section~\ref{sec:exp_mbb_beam} tests robustness across mesh resolutions; Section 5.2 explores filtering effects; Section 5.3 introduces algorithm switching and an updated MMA optimizer; Section 5.4 extends to 3D TO; Section 5.5 quantifies matrix assembly gains; Section 5.6 highlights AD in sensitivity analysis; and Section 5.7 assesses multi-backend and graphics processing units (GPUs) acceleration benefits. These examples underscore SOPTX's strengths and its potential in research and engineering.

\subsection{MBB Beam}\label{sec:exp_mbb_beam}
\begin{figure}[htb]
	\centering
	\includegraphics[width=0.6\textwidth]{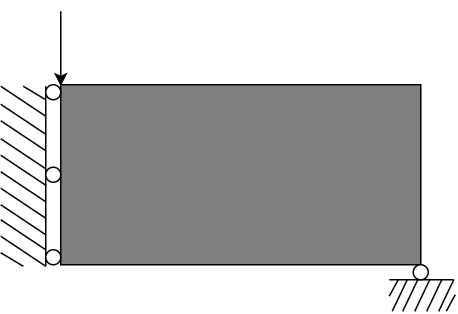}
	\caption{Geometry of the MBB beam: hinged at the left edge and bottom right corner, with a downward concentrated load applied at the top left.}
	\label{fig:mbb_beam}
\end{figure}

As introduced, we first validate SOPTX using the MBB beam, a classical benchmark in TO widely used to assess optimization algorithms. This section minimizes the structural compliance of the MBB beam (see Figure~\ref{fig:mbb_beam}), demonstrating SOPTX's adaptability to various PDE models and meshes. The beam is hinged at the left edge and bottom right corner, with horizontal displacements constrained, and a downward load $T = -1$ applied at the upper left corner. The target volume fraction is 0.5, with material properties $E = 1$ and $\nu = 0.3$, penalization factor $p = 3$, and filter radius $r = 6.0$ to ensure smooth topology and suppress checkerboard patterns.

Compared to the cantilever beam in Section~\ref{sec:exp_canti_beam}, SOPTX enables a seamless switch to the MBB beam problem. Users can modify the PDE model as follows:
\begin{lstlisting}[language=Python]
from soptx.pde import MBBBeam2dData1
pde = MBBBeam2dData1(xmin=0, xmax=150, ymin=0, ymax=50, T=-1)
\end{lstlisting}
The complete model definition is provided in Appendix~\ref{sec:code_mbb}.

\begin{figure}[htb]
	\centering
	\subfloat[$c(\rho)=219.5199,~v(\rho)=0.5001$]{
		\includegraphics[width=0.5\textwidth]{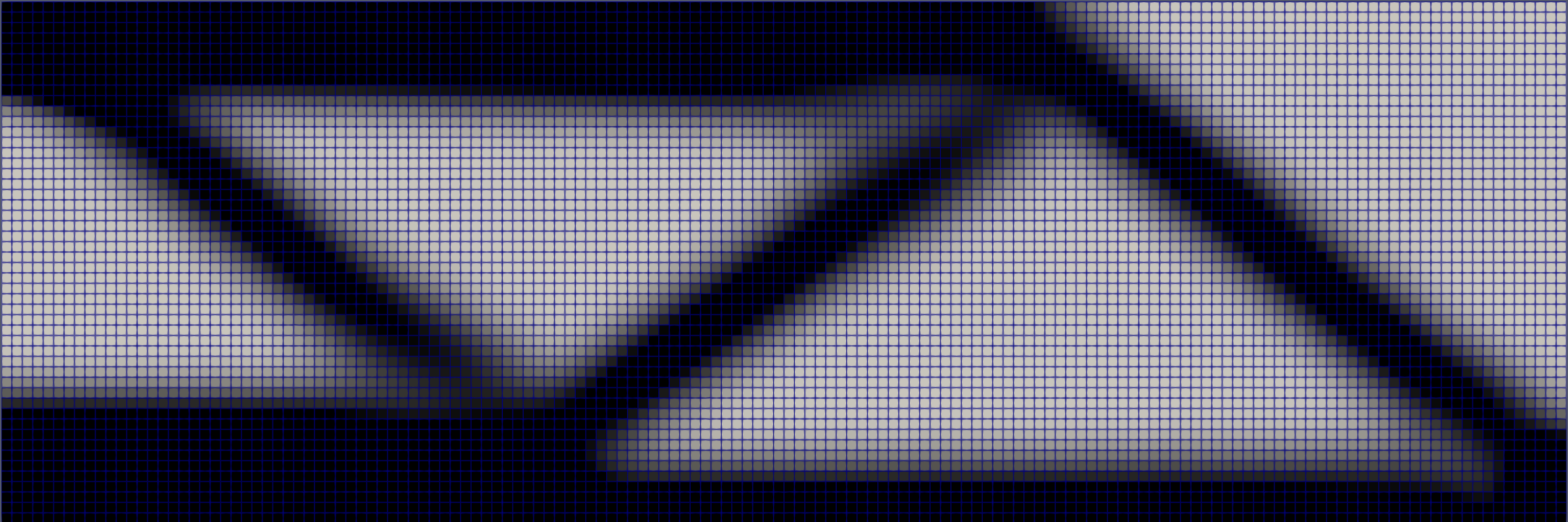}
		\label{fig:mbb_2d_u2}
	}
	\subfloat[$c(\rho)=217.6605,~v(\rho)=0.4999$]{
		\includegraphics[width=0.5\textwidth]{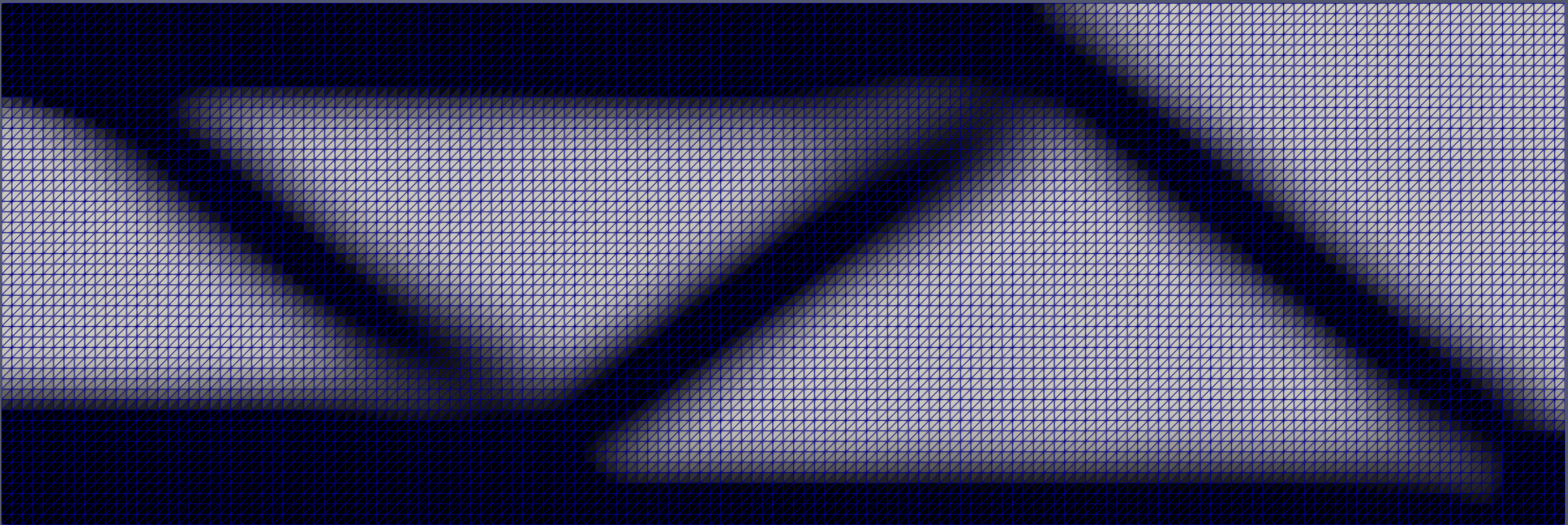}
		\label{fig:mbb_2d_tri}
	}
	\caption{Optimized topologies of the MBB beam using a uniform quadrilateral mesh (left) and a triangular mesh (right).}
	\label{fig:mbb_2d_all}
\end{figure}

To demonstrate SOPTX's adaptability across mesh types, we perform TO on a $150 \times 50$ uniform quadrilateral mesh and a triangular mesh, both initialized with a uniform material density equal to the target volume fraction 0.5. The triangular mesh is generated using:
\begin{lstlisting}[language=Python]
from fealpy.mesh import TriangleMesh
mesh = TriangleMesh.from_box([0, 150, 0, 50], nx=150, ny=50)
\end{lstlisting}

The optimized topologies, shown in Figure~\ref{fig:mbb_2d_all}, are consistent across both mesh types, with compliance and volume fraction differences below $1\%$. This consistency underscores SOPTX's mesh-independence and algorithmic robustness, facilitated by its modular design, which allows easy switching between PDE models and mesh configurations.

\subsection{Different Filtering Methods}
In TO, filtering methods are essential for smoothing design variables, controlling structural details, and ensuring the quality and manufacturability of the optimized results. Thanks to its modular architecture, SOPTX allows users to seamlessly switch between different filtering strategies by simply replacing the filter class, without modifying other components of the framework.

This section compares the results of two common filters (the density filter and the Heaviside projection filter) applied to the MBB beam problem from Section~\ref{sec:exp_mbb_beam}. All parameters are kept identical except for the filter choice, with the default filter radius set to $r=6.0$.

\begin{figure}[htb]
	\centering
	\subfloat[$c(\rho)=235.7337,~v(\rho)=0.5000$]{
		\includegraphics[width=0.5\textwidth]{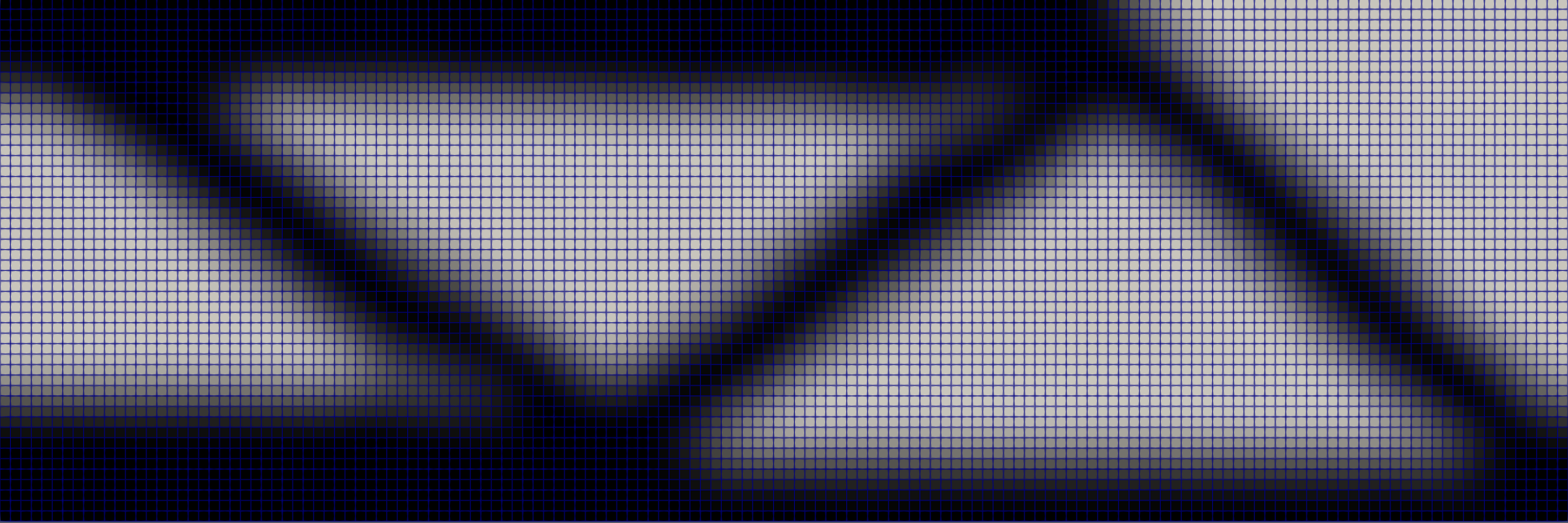}
		\label{fig:mbb_dens}
	}
	\subfloat[$c(\rho)=191.4873,~v(\rho)=0.5000$]{
		\includegraphics[width=0.5\textwidth]{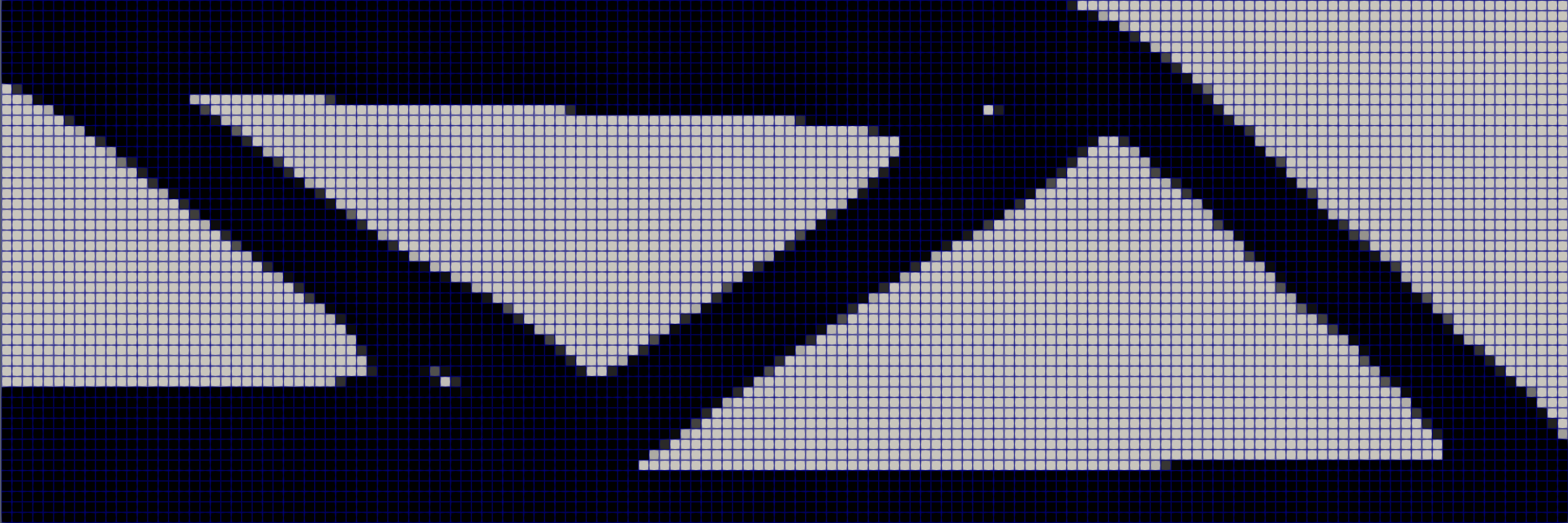}
		\label{fig:mbb_heavi}
	}
	\caption{Optimized topologies of the MBB beam using the density filter (left) and the Heaviside projection filter (right).}
	\label{fig:mbb_filter_all}
\end{figure}

The density filter smooths the design variables through weighted averaging, eliminating small-scale features and producing a gradual material transition. This is ideal for designs requiring structural continuity. In SOPTX, it is applied via:
\begin{lstlisting}[language=python]
dens_filter = DensityBasicFilter(mesh=mesh, rmin=6.0)
\end{lstlisting}

In contrast, the Heaviside projection filter builds on the density filter by adding a projection step. It gradually increases the projection parameter $\beta$ to drive variables toward 0 or 1, yielding a clear black-and-white topology suited for strict manufacturability. To avoid overly thick structures early on, the filter radius is reduced to $r=4.5$. It is enabled in SOPTX with:
\begin{lstlisting}[language=python]
heavi_filter = HeavisideProjectionBasicFilter(mesh=mesh, rmin=4.5, beta=1, max_beta=512, continuation_iter=50)   
\end{lstlisting}

The optimized topologies are shown in Figure~\ref{fig:mbb_filter_all}. The density filter yields a smoother design with blurred edges, while the Heaviside projection filter produces a crisp, binarized layout with lower compliance (191.4873 vs. 235.7337), indicating a stiffer structure. However, with the Heaviside projection filter, small holes may still appear due to its focus on global binarization rather than strict local control.

In summary, SOPTX's flexibility allows efficient exploration of filter impacts on design smoothness, manufacturability, and performance.

\subsection{Different Optimization Algorithms}
In TO, the choice of optimization algorithm impacts both computational efficiency and design quality. SOPTX's modular design allows seamless switching between optimizers, such as from the Optimality Criteria (OC) method to the more versatile Method of Moving Asymptotes (MMA). This section demonstrates this process using the MBB beam problem and highlights the refactored MMA in SOPTX.

\begin{figure}[htb]
	\centering
	\subfloat[$c(\rho)=219.4545,~v(\rho)=0.5000$]{
		\includegraphics[width=0.5\textwidth]{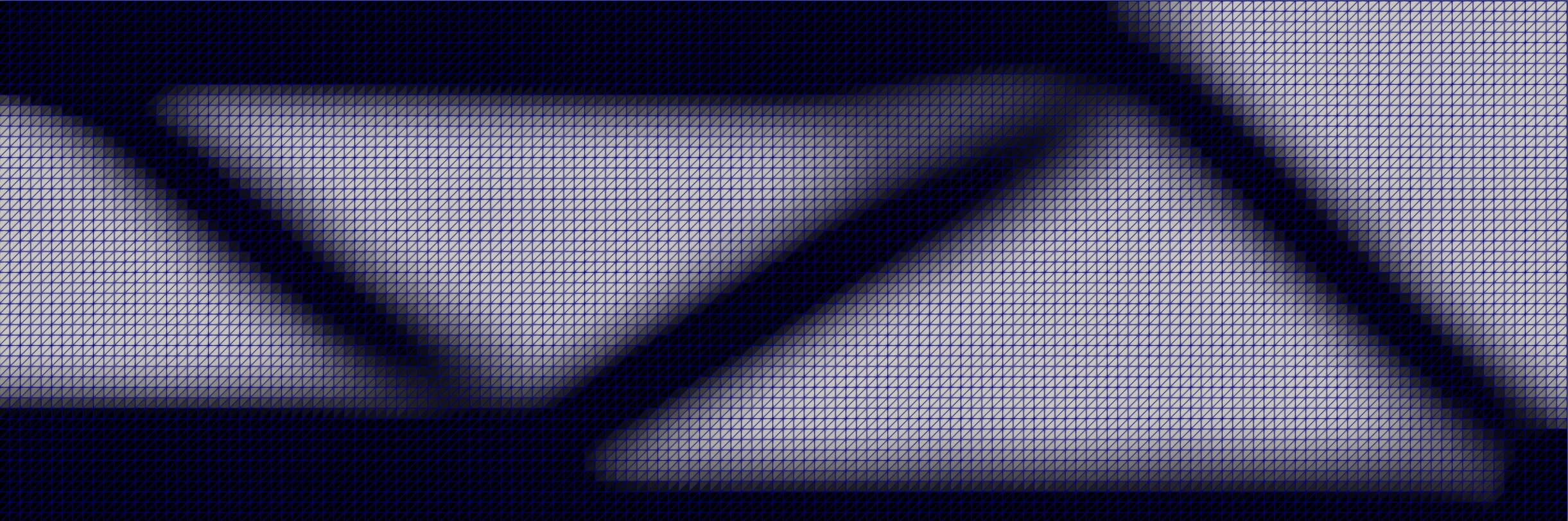}
		\label{fig:mbb_mma_u2}
	}
	\subfloat[$c(\rho)=220.4377,~v(\rho)=0.4995$]{
		\includegraphics[width=0.5\textwidth]{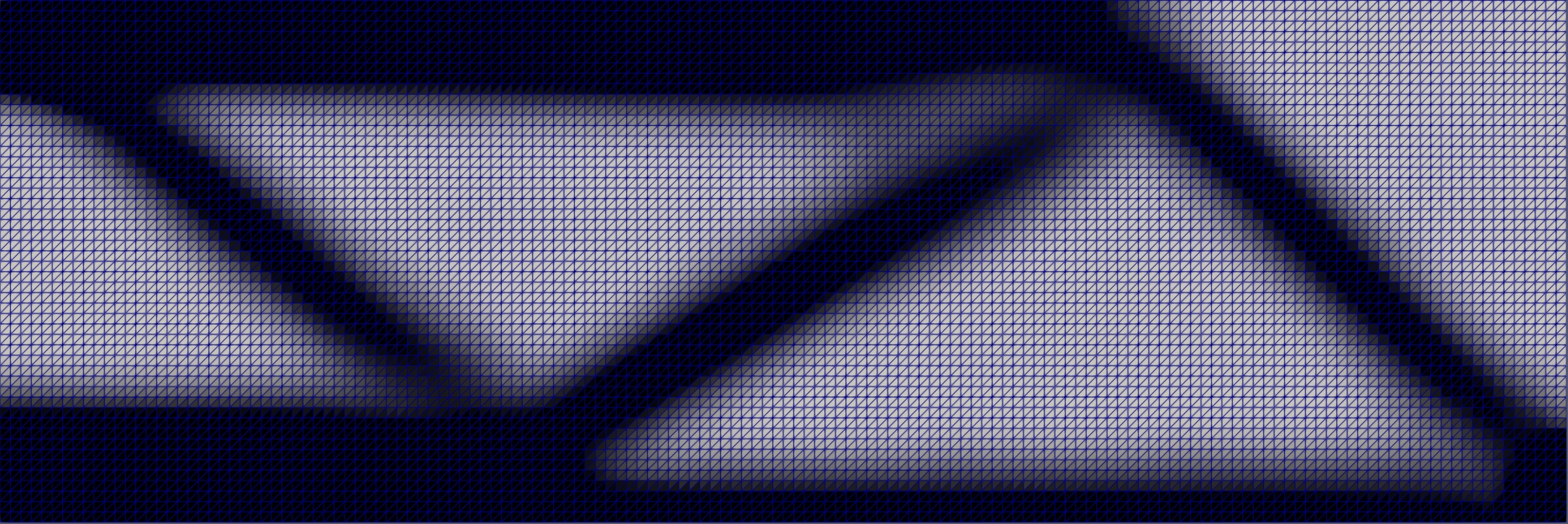}
		\label{fig:mbb_mma_tri}
	}
	\caption{Optimized topologies of the MBB beam using MMA optimizer on a structured quadrilateral mesh (left) and a triangular mesh (right).}
	\label{fig:mbb_mma_all}
\end{figure}

To switch to MMA, users simply import and configure the \texttt{MMAOptimizer} class:
\begin{lstlisting}[language=python]
from soptx.opt import MMAOptimizer
optimizer = MMAOptimizer(objective=objective, 
		constraint=constraint, filter_=sens_filter,
		options={'max_iterations': 200, 'tolerance': 0.01})
\end{lstlisting}

For the MBB beam problem, MMA produces topologies nearly identical to OC's, with compliance and volume fraction differences below $1\%$, as shown in Figure~\ref{fig:mbb_mma_all}. This consistency confirms MMA's reliability for single-constraint problems.

SOPTX's MMA, refactored from Krister Svanberg's implementation ~\cite{Svanberg2007MmaAG}, allows users to adjust internal parameters, including \texttt{m}, \texttt{n}, \texttt{xmin}, \texttt{xmax}, and control parameters like \texttt{a0}, \texttt{a}, \texttt{c}, \texttt{d}. This provides greater flexibility than traditional "black-box" versions, enhancing performance and adaptability across platforms. An example configuration is shown below:
\begin{lstlisting}[language=python]
optimizer.options.set_advanced_options(m=1, n=NC,
			xmin=bm.zeros((NC, 1)), xmax=bm.ones((NC, 1)),
			a0=1, a=bm.zeros((1, 1)),
			c=1e4 * bm.ones((1, 1)), d=bm.zeros((1, 1)))
\end{lstlisting}

Through this refactoring, SOPTX not only enhances the performance of the MMA algorithm but also overcomes the limitations of conventional implementations. It provides users with greater control and cross-platform flexibility, making it particularly advantageous for both research and engineering applications in TO.

\subsection{3D Extension}\label{sec:exp_canti3d}
\begin{figure}[htbp]
	\centering
	\includegraphics[width=0.5\textwidth]{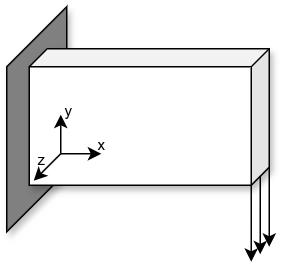}
	\caption{Geometric configuration of the 3D cantilever beam. The left end is fully fixed, and a downward concentrated load is applied at the bottom of the right end.}
	\label{fig:cantilever_3d}
\end{figure}

SOPTX's modular design enables a smooth transition from 2D to 3D TO. This section showcases its capability with a 3D cantilever beam compliance minimization problem (see Figure~\ref{fig:cantilever_3d}). The beam is fixed on the left and bears a downward load $T=-1$ at the bottom right. A $60\times20\times4$ hexahedral mesh is applied, with a target volume fraction of 0.3, material properties $E=1$ and $\nu=0.3$, and penalization factor $p=3$. A sensitivity filter (radius $r=1.5$) ensures smooth results.

Compared to the 2D cantilever beam in Section~\ref{sec:exp_canti_beam}, SOPTX extends to 3D by simply replacing the PDE model and mesh configuration. The 3D cantilever beam model is implemented as follows:
\begin{lstlisting}[language=python]
from soptx.pde import Cantilever3dData1
pde = Cantilever3dData1(xmin=0, xmax=60, ymin=0, ymax=20, zmin=0, zmax=4, T=-1)
mesh = UniformMesh(extent=[0, 60, 0, 20, 0, 4], h=[1, 1, 1], origin=[0, 0, 0])
\end{lstlisting}
The complete definition of the \texttt{Cantilever3dData1} model is available in Appendix~\ref{sec:code_canti_3d}.

\begin{figure}[htbp]
	\centering
	\includegraphics[width=1.0\textwidth]{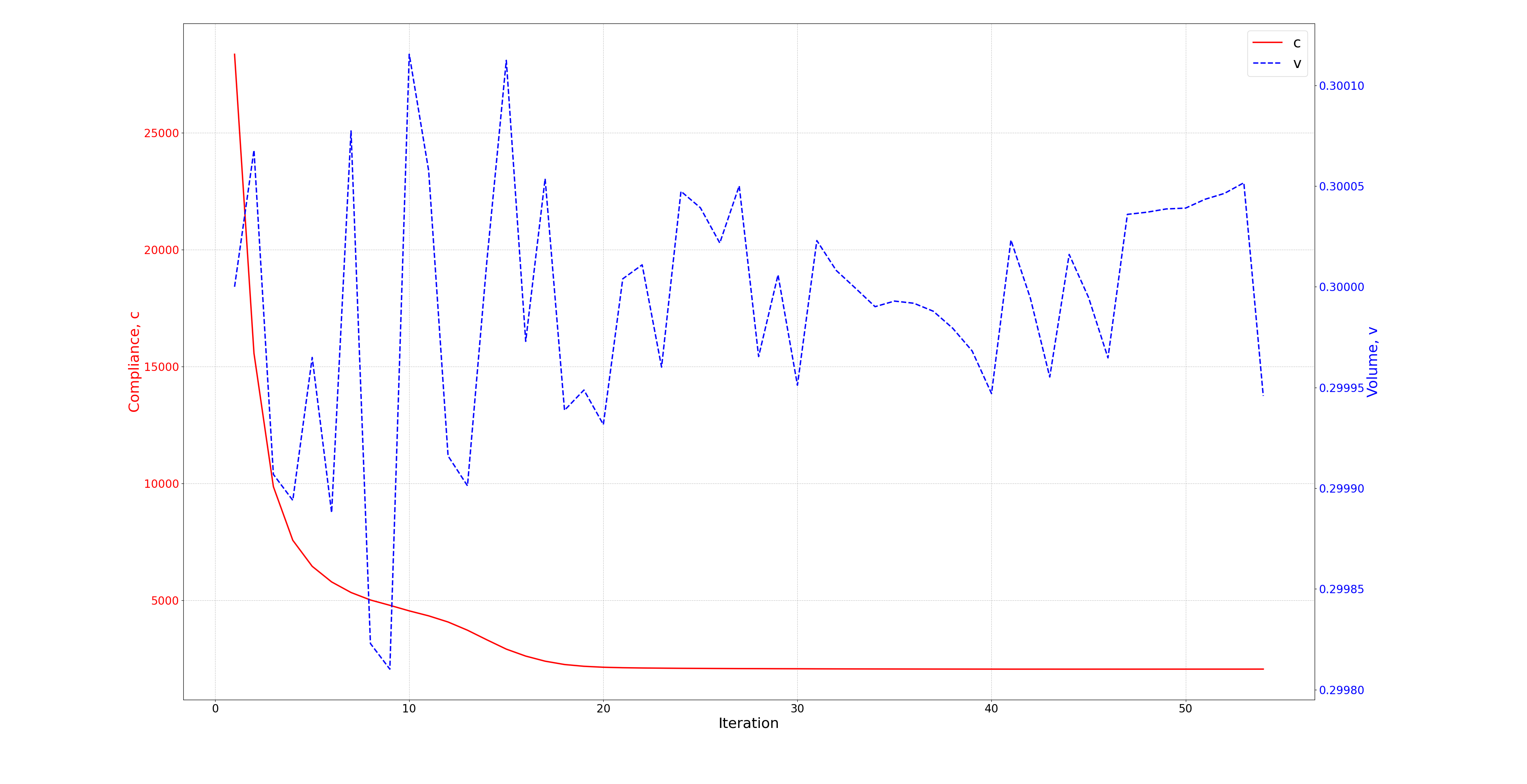}
	\caption{Convergence histories of the compliance $c(\rho)$ and volume fraction $v(\rho)$ for the 3D cantilever beam initialized with a uniform density of 0.3.}
	\label{fig:canti_3d_convergence}
\end{figure}

The initial material density is uniformly set to the target volume fraction of 0.3. Using the OC optimizer and a sensitivity filter with radius $r=1.5$, the compliance $c(\rho)$ shows a rapid initial decrease, followed by a gradual convergence with slight fluctuations, stabilizing at around 2,000 by iteration 54 (see Figure~\ref{fig:canti_3d_convergence}). The volume fraction remains nearly constant around 0.3 throughout the process, indicating effective constraint control.

\begin{figure}[htb]
	\centering
	\subfloat[$c(\rho)=5338.9462,~v(\rho)=0.3001$]{
		\includegraphics[width=0.33\textwidth]{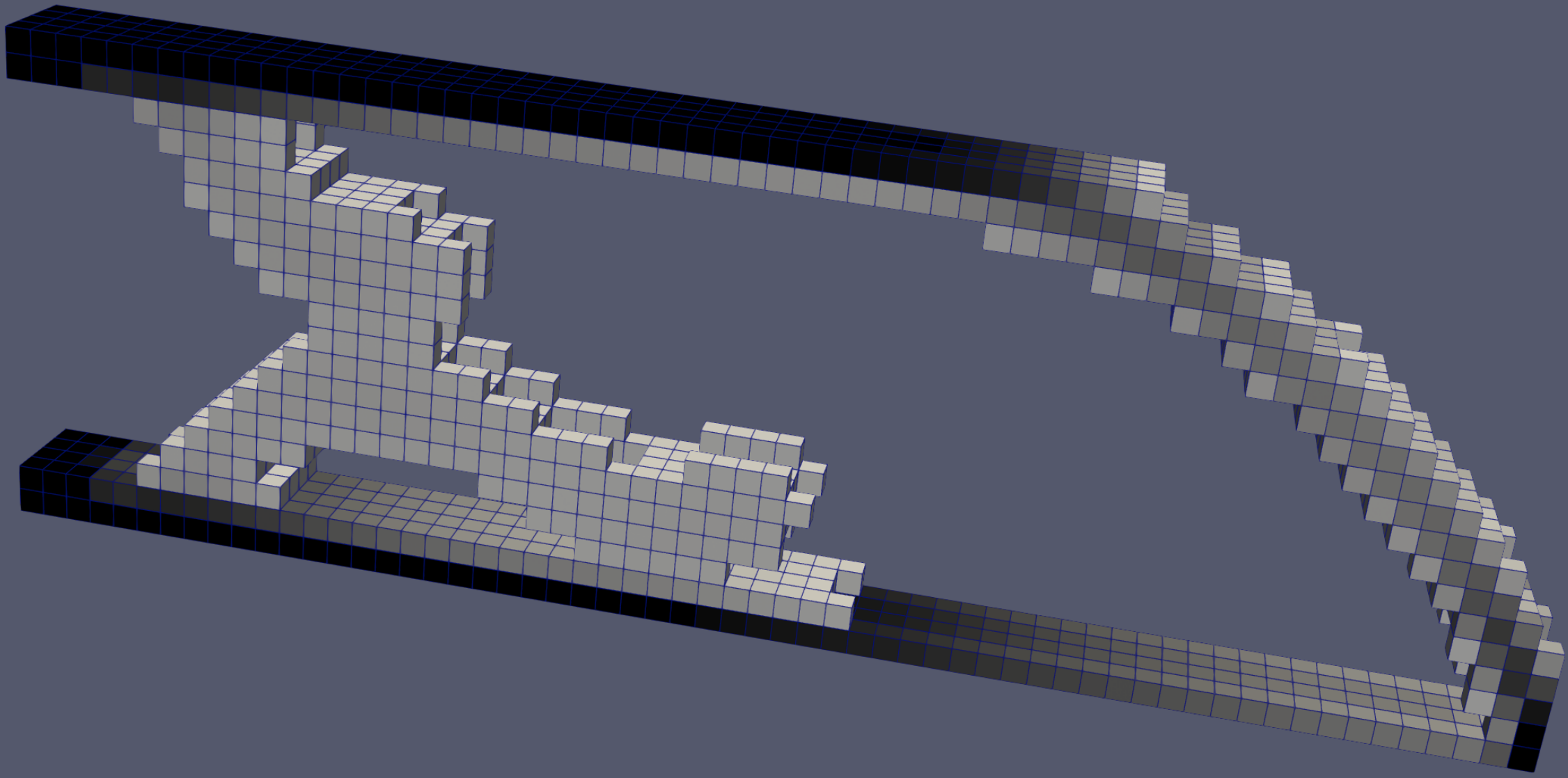}
		\label{fig:canti3d_1}
	}
	\subfloat[$c(\rho)=2122.6707,~v(\rho)=0.3000$]{
		\includegraphics[width=0.33\textwidth]{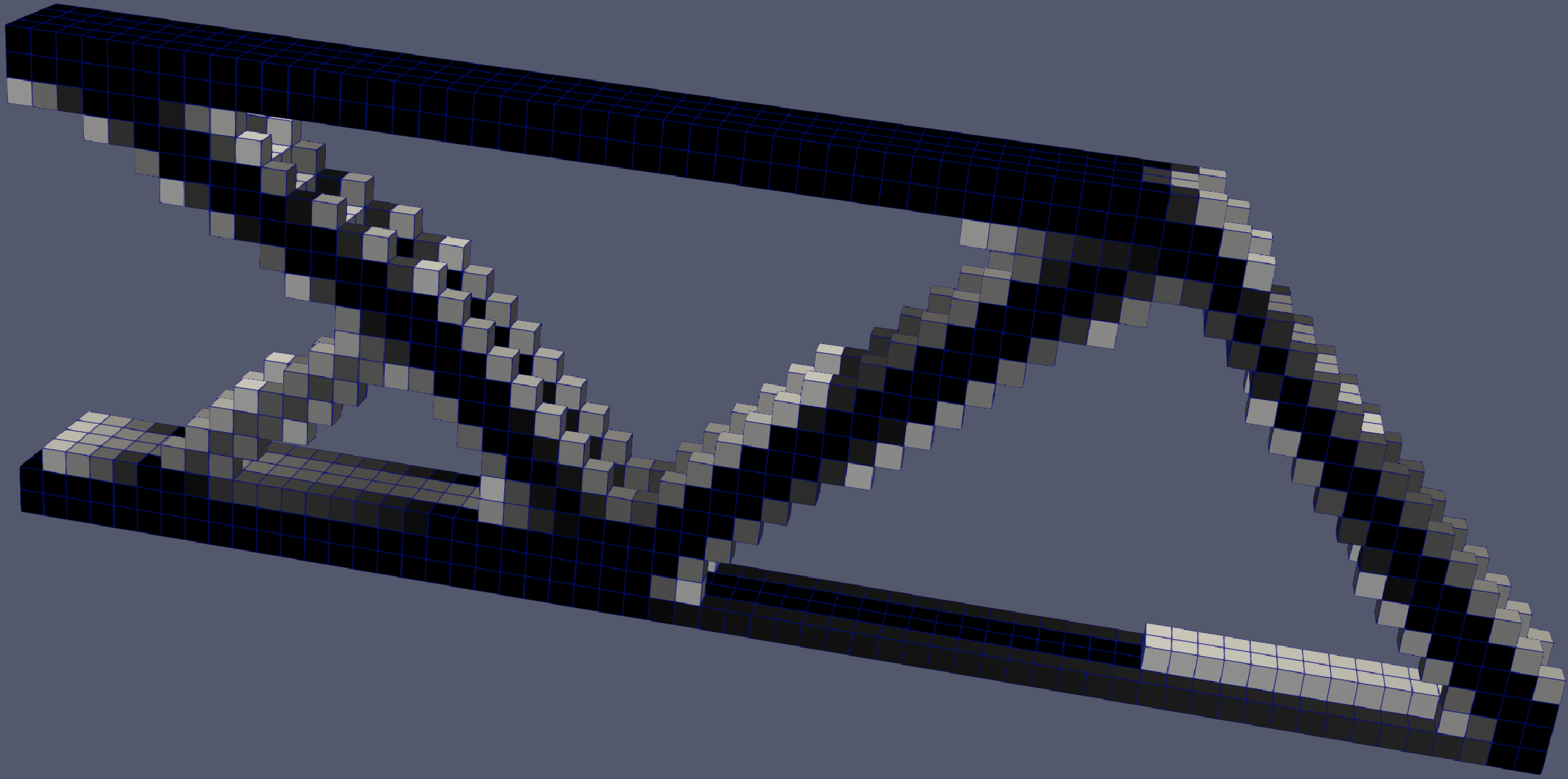}
		\label{fig:canti3d_2}
	}
	\subfloat[$c(\rho)=2063.5625,~v(\rho)=0.2999$]{
		\includegraphics[width=0.33\textwidth]{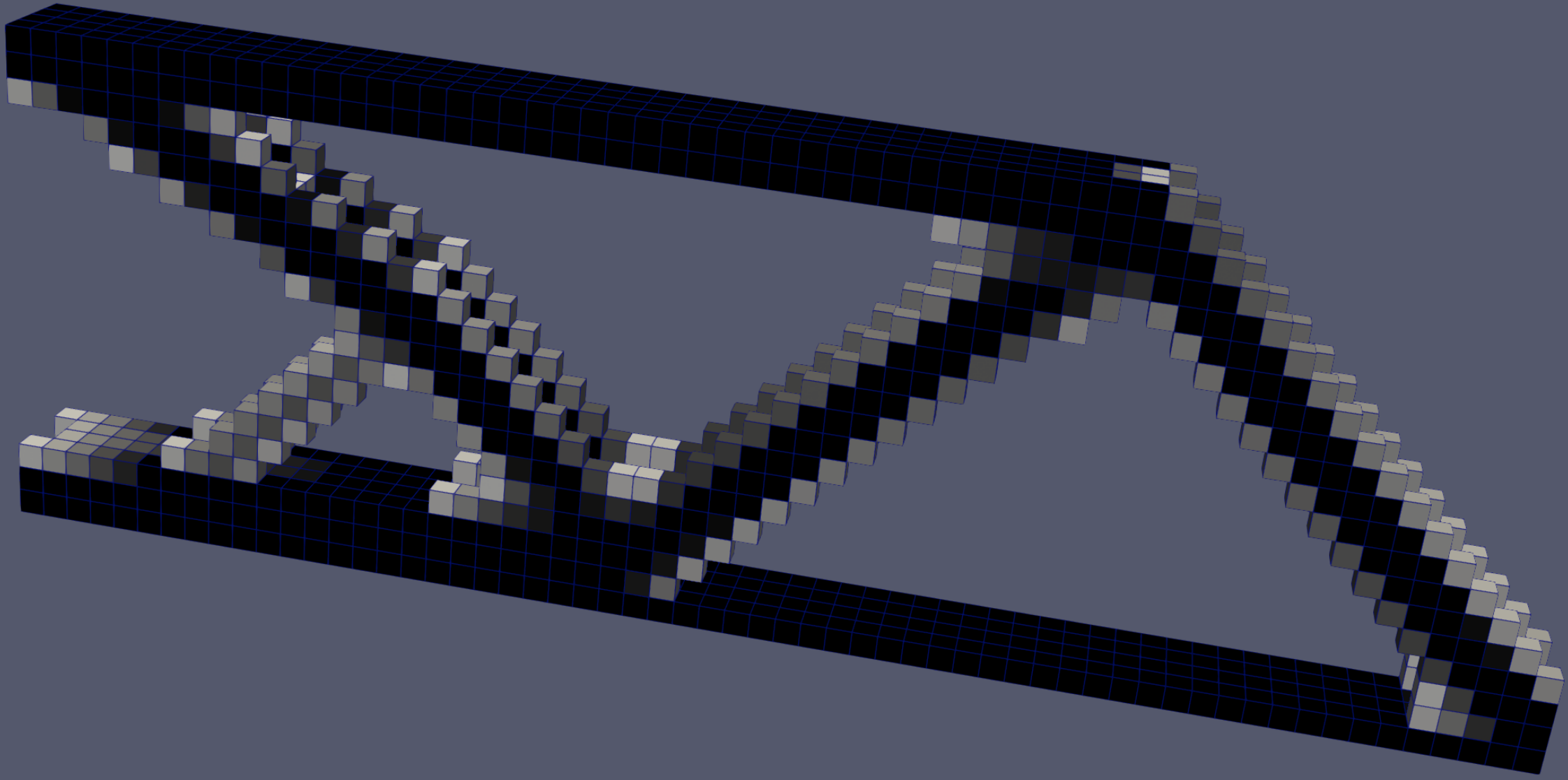}
		\label{fig:canti3d_3}
	}
	\caption{Topology layouts at iterations 7 (left), 21 (middle), and 54 (right) during the optimization of the 3D cantilever beam. Elements with ($\rho > 0.3$) are visualized. Each subfigure also reports the compliance and volume fraction.}
	\label{fig:canti3d_all}
\end{figure}

To better visualize the optimization results, Figure~\ref{fig:canti3d_all} shows the topology layouts at iterations 7, 21, and 54. Only elements with density values $\rho>0.3$ are rendered to highlight the structural features.

3D TO increases computational complexity significantly. SOPTX addresses this with efficient matrix assembly and multi-backend support, enhancing performance for large-scale problems. Details follow in subsequent sections

\subsection{Application of Fast Matrix Assembly}
In 3D TO problems, computational efficiency is critical. For the 3D cantilever beam example (Section~\ref{sec:exp_canti3d}), the mesh of $60\times20\times4$ elements leads to 19,215 displacement degrees of freedom and 4,800 density variables, requiring 54 iterations. Traditional matrix assembly in finite element methods involves redundant computations, creating a bottleneck.

SOPTX addresses this with a fast matrix assembly technique that separates element-dependent and element-independent parts, reusing invariant quantities to accelerate the process. Additionally, symbolic integration precomputation analytically integrates invariant terms in advance, improving accuracy and stability.
\begin{itemize} 
	\item \textbf{Fast assembly}: Accelerates matrix assembly by separating element-dependent and element-independent terms, coupled with efficient numerical integration.
	\item \textbf{Symbolic Fast Assembly}: Further improves accuracy and stability by analytically integrating invariant terms ahead of time.
\end{itemize}

Users can easily switch between assembly methods by setting the \texttt{assembly\_method} parameter in the solver initialization, as shown below:
\begin{lstlisting}[language=python]
# Fast Assembly
solver = ElasticFEMSolver(..., assembly_method=AssemblyMethod.FAST, ...)
# Symbolic Fast Assembly
solver = ElasticFEMSolver(..., assembly_method=AssemblyMethod.SYMBOLIC, ...)
\end{lstlisting}

To evaluate the efficiency of different matrix assembly techniques in SOPTX, we use the same 3D cantilever beam optimization problem from Section~\ref{sec:exp_canti3d} as a representative test case. Table~\ref{tab:assembly_comparison} presents a performance comparison of these techniques
\begin{table}[htbp]
	\centering
	\caption{Performance comparison of matrix assembly techniques in the 3D cantilever beam optimization problem. All values are reported in seconds.}
	\begin{tabular}{cccccc}
		\toprule
		\textbf{Assembly Technique} & \textbf{Total} & \textbf{1st Iter.} & \textbf{1st Assembly} & \textbf{Avg. Iter.} & \textbf{Avg. Assembly} \\
		\midrule
		Original Assembly & 68.605 & 3.342 & 1.197 & 1.231 & 0.838 \\
		Fast Assembly & 39.826 & 2.387 & 0.342 & 0.706 & 0.276 \\
		Symbolic Fast Assembly & 41.194 & 5.877 & 3.853 & 0.666 & 0.272 \\
		\bottomrule
	\end{tabular}
	\label{tab:assembly_comparison}
\end{table}

As shown in Table~\ref{tab:assembly_comparison}, the original assembly method exhibits an average assembly time of \SI{0.838}{s} per iteration, accounting for $68\%$ of the total iteration time and thus creating a significant bottleneck. The fast assembly technique substantially improves this, reducing the average assembly time to \SI{0.276}{s}, which is only $33\%$ of the original assembly time. The symbolic fast assembly technique, while incurring a longer first iteration time (\SI{3.853}{s} for the first assembly), achieves the fastest average assembly time of \SI{0.272}{s} with enhanced accuracy. This performance difference arises because symbolic fast assembly requires precomputing symbolic integrals during the first iteration, which increases the initial iteration time but allows the reuse of intermediate results in subsequent iterations, thereby significantly reducing assembly time in later stages.

In summary, SOPTX's fast assembly strategies significantly reduce computational cost and improve accuracy, making it ideal for large-scale TO problems.

\subsection{Application of Automatic Differentiation}\label{sec:exp_canti3d_ad}
Sensitivity analysis is vital in TO to update design variables. Traditional methods rely on manual sensitivity derivation, which is time-consuming and error-prone, especially for complex models. The SOPTX framework introduces AD to streamline this process, enabling users to prioritize problem modeling over mathematical derivations. Here, we demonstrate AD's application using the 3D cantilever beam example from Section~\ref{sec:exp_canti3d}.

The SOPTX framework supports multiple computational backends(e.g., NumPy, PyTorch, JAX). To enable AD, users must switch from the default NumPy backend to an AD-enabled backend like PyTorch~\cite{paszke2017automatic} with a single line of code:
\begin{lstlisting}[language=python]
bm.set_backend('pytorch')
\end{lstlisting}

Sensitivity computation is controlled via the \texttt{diff\_mode} parameter: \texttt{'manual'} for hand-derived formulas and \texttt{'auto'} for AD. For example, AD can be enabled for the compliance objective while keeping manual differentiation for the volume constraint:
\begin{lstlisting}[language=python]
obj_config = ComplianceConfig(diff_mode='auto')
objective = ComplianceObjective(solver=solver, config=obj_config)
cons_config = VolumeConfig(diff_mode='manual')
constraint = VolumeConstraint(solver=solver, volume_fraction=0.5, config=cons_config)
\end{lstlisting}
This setup allows selective use of AD or manual differentiation. Further details on AD implementation are in Section~\ref{sec:ad_compliance}.

\begin{figure}[htb] 
	\centering 
	\includegraphics[width=1.0\textwidth]{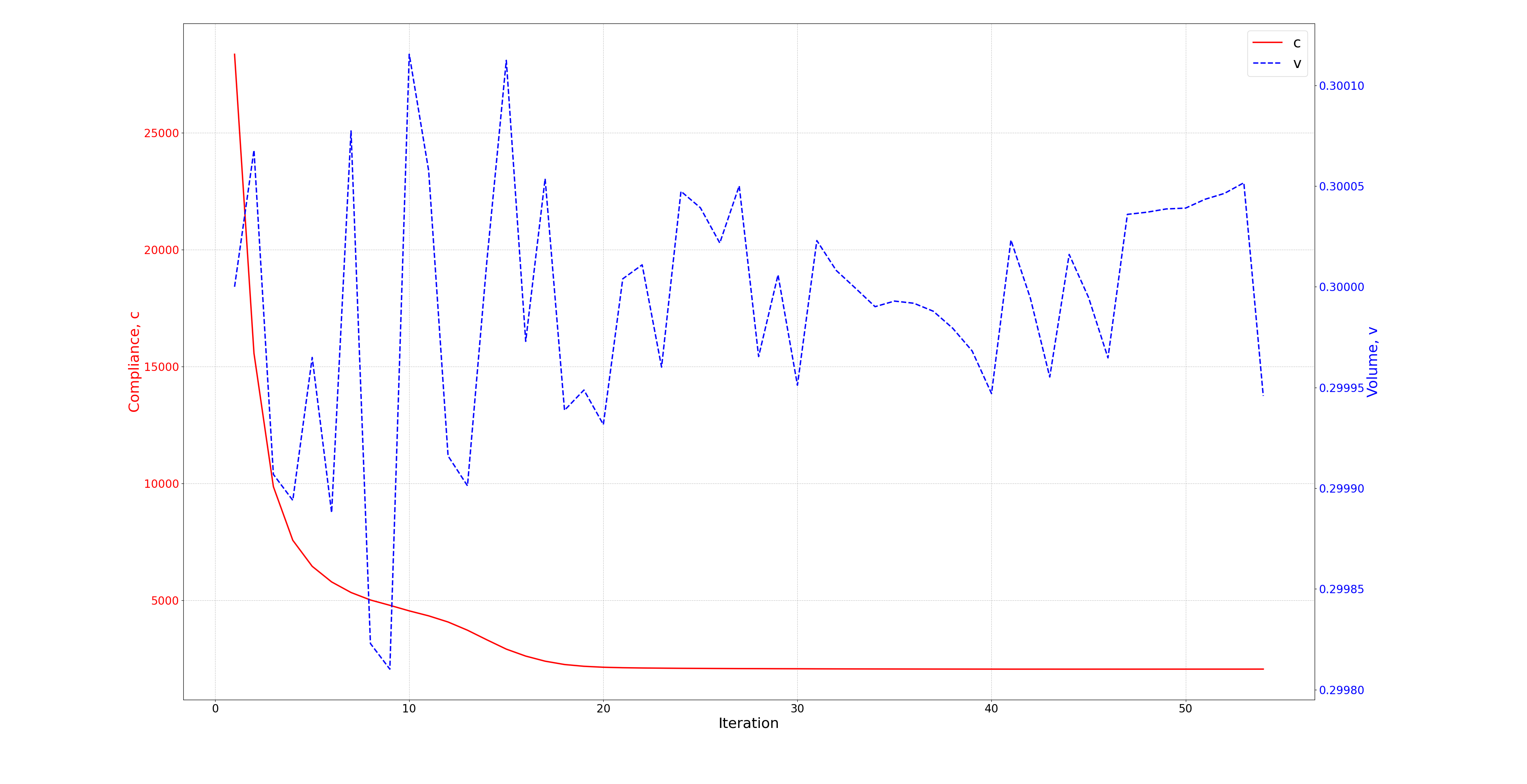} \caption{Convergence histories of the compliance $c(\rho)$ and volume fraction $v(\rho)$ for the 3D cantilever beam optimization using automatic differentiation.} \label{fig:canti3d_ad_convergence} 
\end{figure}

To verify the correctness and effectiveness of AD, we conduct tests using the 3D cantilever beam problem described in Section~\ref{sec:exp_canti3d}. All parameters are kept identical. The structural compliance sensitivity is computed using AD. Figure~\ref{fig:canti3d_ad_convergence} shows the convergence histories of the compliance and volume fraction.

\begin{figure}[htb]
	\centering
	\subfloat[$c(\rho)=5338.9462,~v(\rho)=0.3001$]{
		\includegraphics[width=0.32\textwidth]{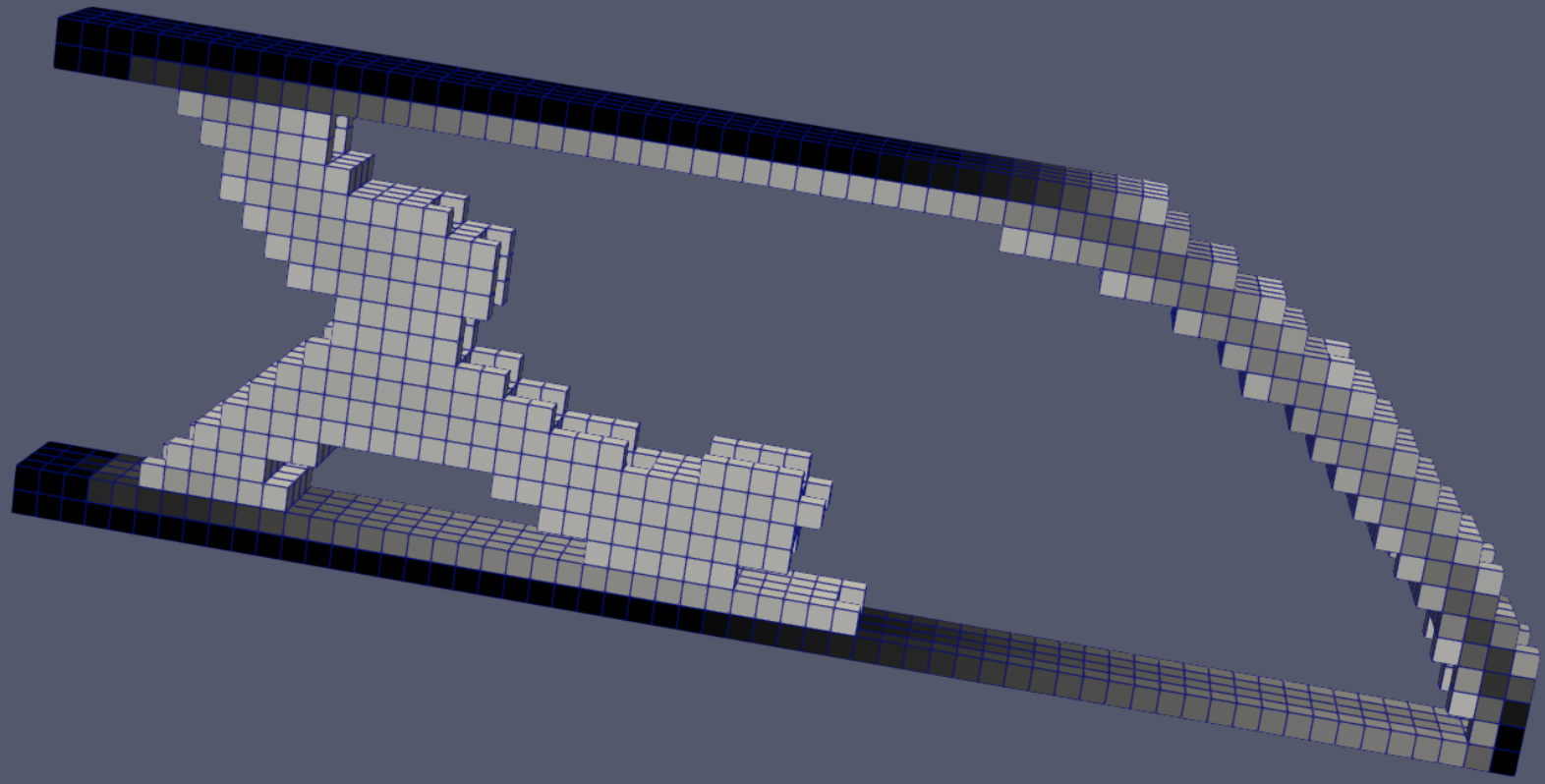}
		\label{fig:canti3d_1_auto}
	}
	\subfloat[$c(\rho)=2122.6707,~v(\rho)=0.3000$]{
		\includegraphics[width=0.32\textwidth]{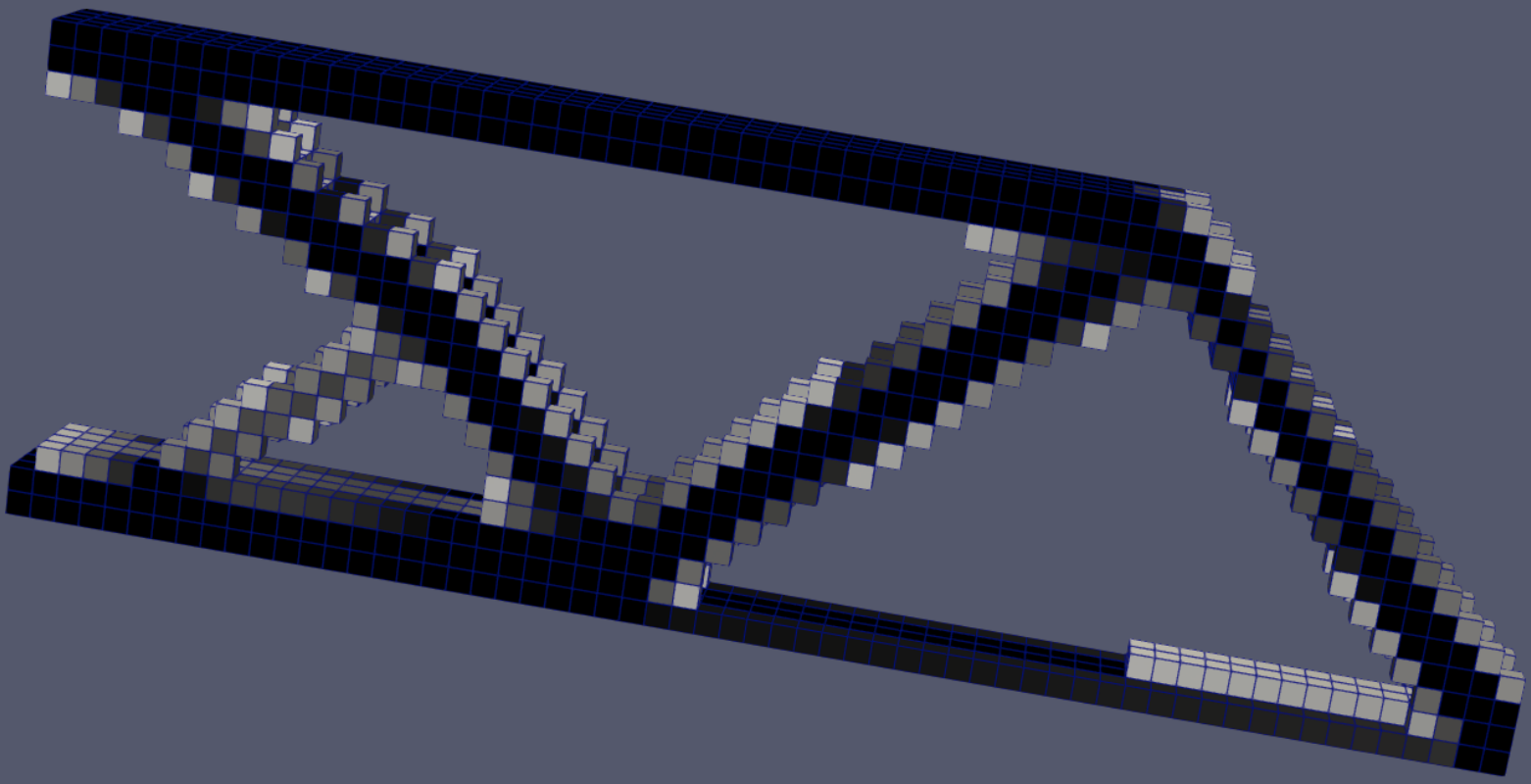}
		\label{fig:canti3d_2_auto}
	}
	\subfloat[$c(\rho)=2063.5625,~v(\rho)=0.2999$]{
		\includegraphics[width=0.32\textwidth]{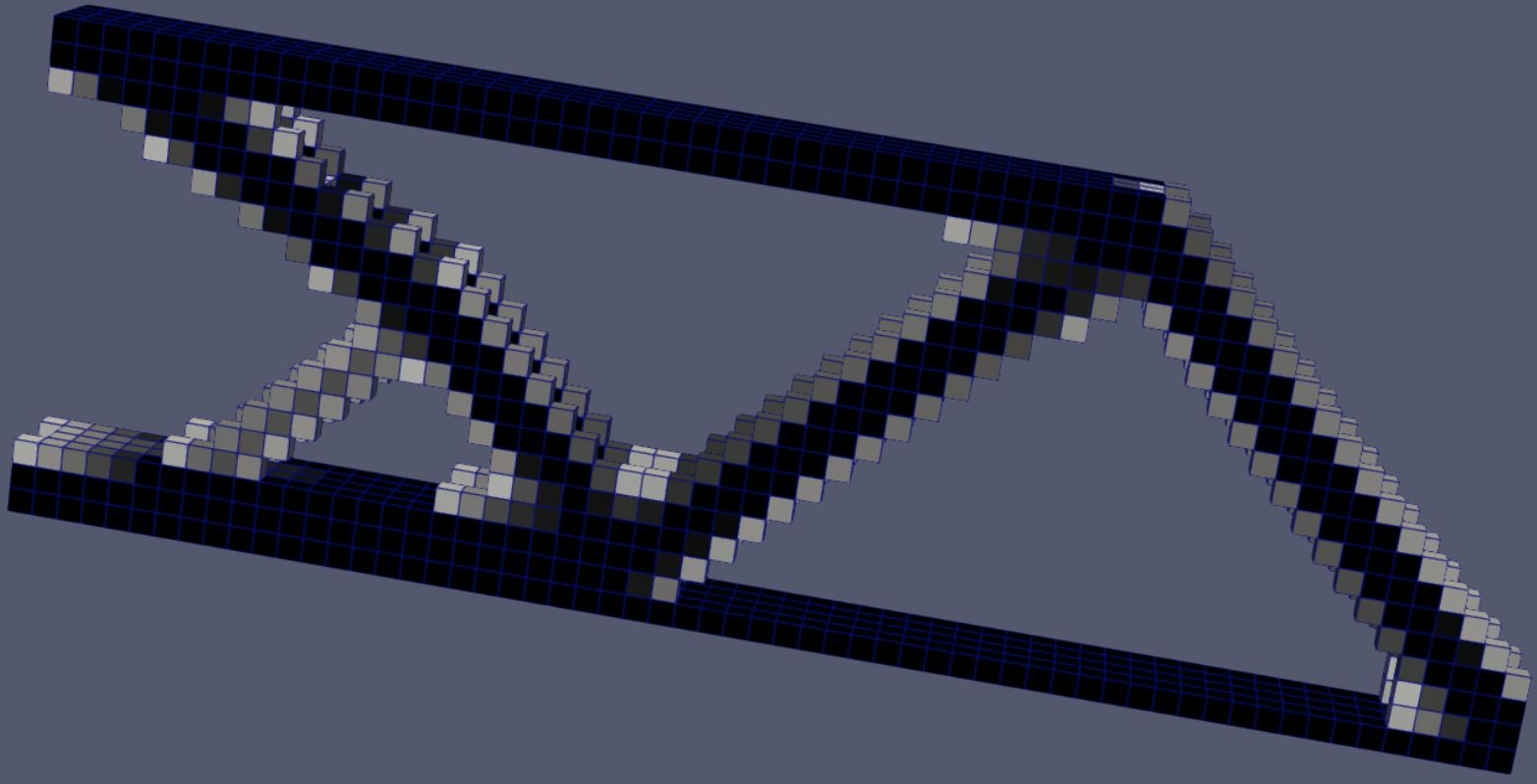}
		\label{fig:canti3d_3_auto}
	}
	\caption{Topology layouts at iterations 7 (left), 21 (middle), and 54 (right) during the 3D cantilever beam optimization using AD. Only elements with $\rho>0.3$ are visualized. Each subfigure also reports the compliance and volume fraction.}
	\label{fig:canti3d_ad_topos}
\end{figure}

To further illustrate the evolution of the topology during optimization, the layouts at selected iterations are shown in Figure~\ref{fig:canti3d_ad_topos}.

As observed from Figure~\ref{fig:canti3d_ad_convergence} and Figure~\ref{fig:canti3d_ad_topos}, the optimization results obtained using AD are fully consistent with those achieved through manual differentiation in Section~\ref{sec:exp_canti3d}. The compliance converges to approximately 2063.5625 after 54 iterations, and the volume fraction remains stably around 0.3. These results confirm that AD accurately computes the sensitivities and effectively guides the optimization process.

To evaluate the computational efficiency of AD, we compare the optimization times of AD and manual differentiation under the same 3D cantilever beam problem setting. As shown in Table~\ref{tab:ad_vs_manual}, the total optimization time and average iteration time for both methods are nearly identical, demonstrating that SOPTX's AD implementation is highly efficient, introducing no significant computational overhead while maintaining the same high precision as manual differentiation.
\begin{table}[htbp]
	\centering
	\setlength{\tabcolsep}{4pt} 
	\caption{Performance comparison of differentiation methods in the 3D cantilever beam optimization problem.}
	\begin{tabular}{ccccc}
		\toprule
		\textbf{Differentiation Method} & \textbf{Iterations} & \textbf{Total Time (s)} & \textbf{1st Iter. Time (s)} & \textbf{Avg. Iter. Time (s)} \\
		\midrule
		Manual Differentiation & 54 & 39.562 & 1.940 & 0.710 \\
		AD & 54 & 39.865 & 1.832 & 0.718 \\
		\bottomrule
	\end{tabular}
	\label{tab:ad_vs_manual}
\end{table}

\begin{figure}[htbp] 
	\centering 
	\includegraphics[width=0.75\textwidth]{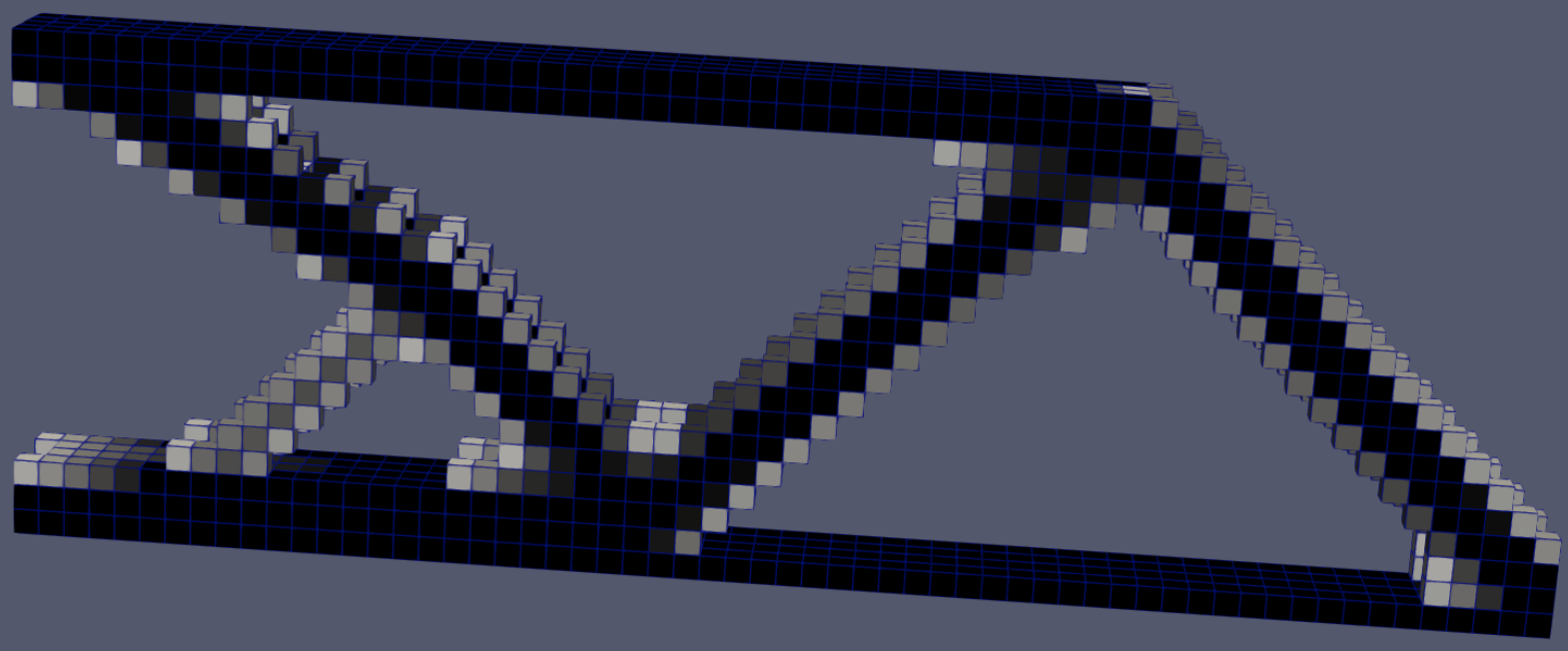} \caption{Final optimized topology of the 3D cantilever beam using AD with the JAX backend.} \label{fig:canti3d_ad_jax} 
\end{figure}

In addition to PyTorch, SOPTX also supports the JAX backend~\cite{jax2018github}, which offers equally powerful AD capabilities. The sensitivity computation code remains identical when using JAX, users simply need to switch to the JAX backend by executing:
\begin{lstlisting}[language=python]
bm.set_backend('jax')
\end{lstlisting}

Using AD with the JAX backend, we perform TO for the 3D cantilever beam, achieving results fully consistent with those from the PyTorch backend. This consistency underscores the stability and flexibility of the SOPTX framework across different computational backends. The final optimized topology is shown in Figure~\ref{fig:canti3d_ad_jax}.

The application of AD in TO offers several significant advantages:
\begin{itemize} 
	\item Simplified model switching: AD automates gradient computations for material interpolation models like SIMP and Rational Approximation of Material Properties (RAMP)~\cite{stolpe2001alternative}, enabling seamless switching without manual derivation and improving development efficiency.
	\item Seamless constraint switching: AD simplifies sensitivity analysis for constraints such as volume, length scale~\cite{guest2009imposing}, connectivity~\cite{li2016structural}, overhangs~\cite{qian2017undercut}, and material usage~\cite{sanders2018multi}, enhancing optimization adaptability.
	\item Support for complex problems: AD efficiently computes gradients in multi-physics or geometrically constrained TO, automatically handling sensitivity computations for techniques like density filtering and Heaviside projection, allowing users to focus on modeling rather than derivations.
\end{itemize}

This section highlights AD's potential in SOPTX, delivering precision and efficiency comparable to manual differentiation without added overhead. SOPTX's multi-backend support enhances its flexibility, while AD simplifies sensitivity analysis and enables exploration of new models and constraints, making SOPTX a valuable tool for topology optimization in education, research, and engineering.

\subsection{Multi-Backend Switching}\label{sec:exp_multi_backend}
SOPTX supports multiple computational backends, including NumPy, PyTorch, and JAX. Users can flexibly switch between them using the \texttt{set\_backend} function.
\begin{lstlisting}[language=python]
bm.set_backend('numpy')    
bm.set_backend('pytorch') 
bm.set_backend('jax')     
\end{lstlisting}

\begin{figure}[htbp]
	\centering
	\subfloat[$c(\rho)=2063.5625,~v(\rho)=0.2999$]{
		\includegraphics[width=0.32\textwidth]{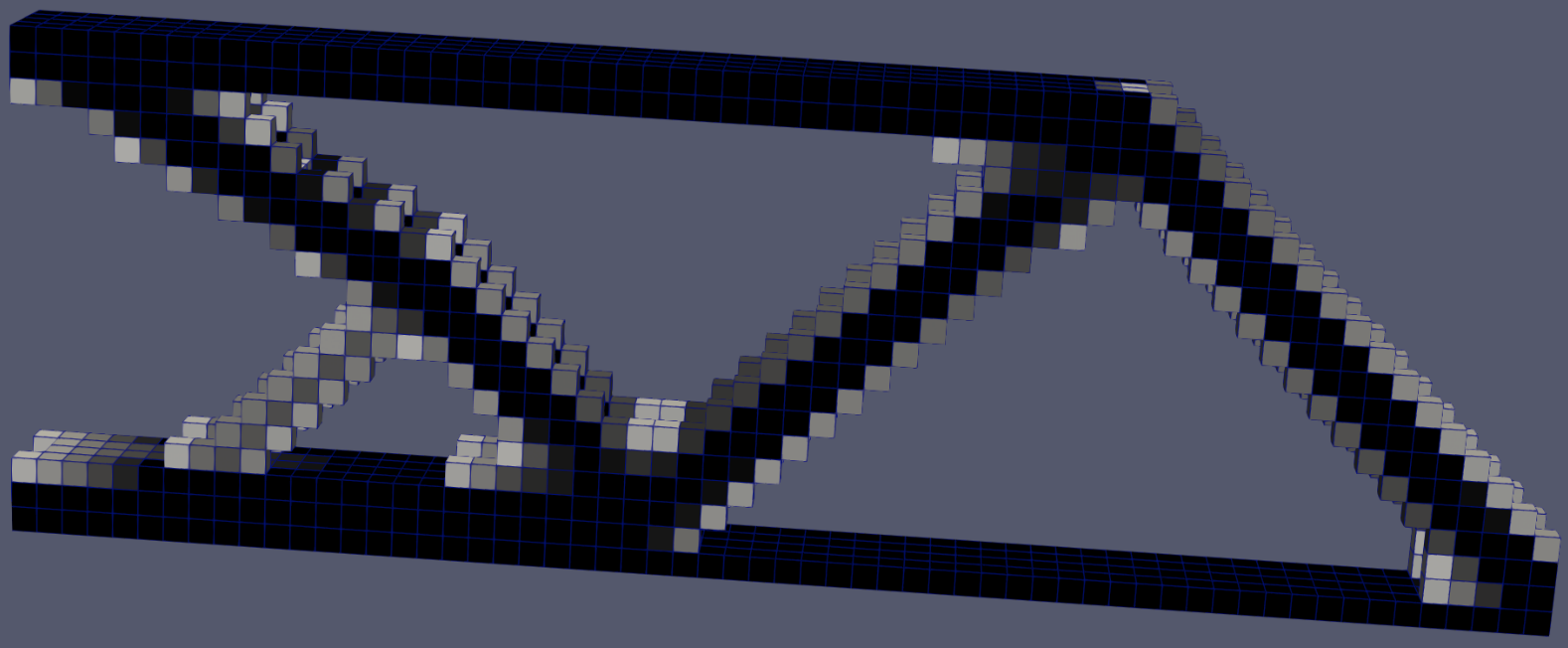}
		\label{fig:canti3d_np}
	}
	\subfloat[$c(\rho)=2063.5625,~v(\rho)=0.2999$]{
		\includegraphics[width=0.32\textwidth]{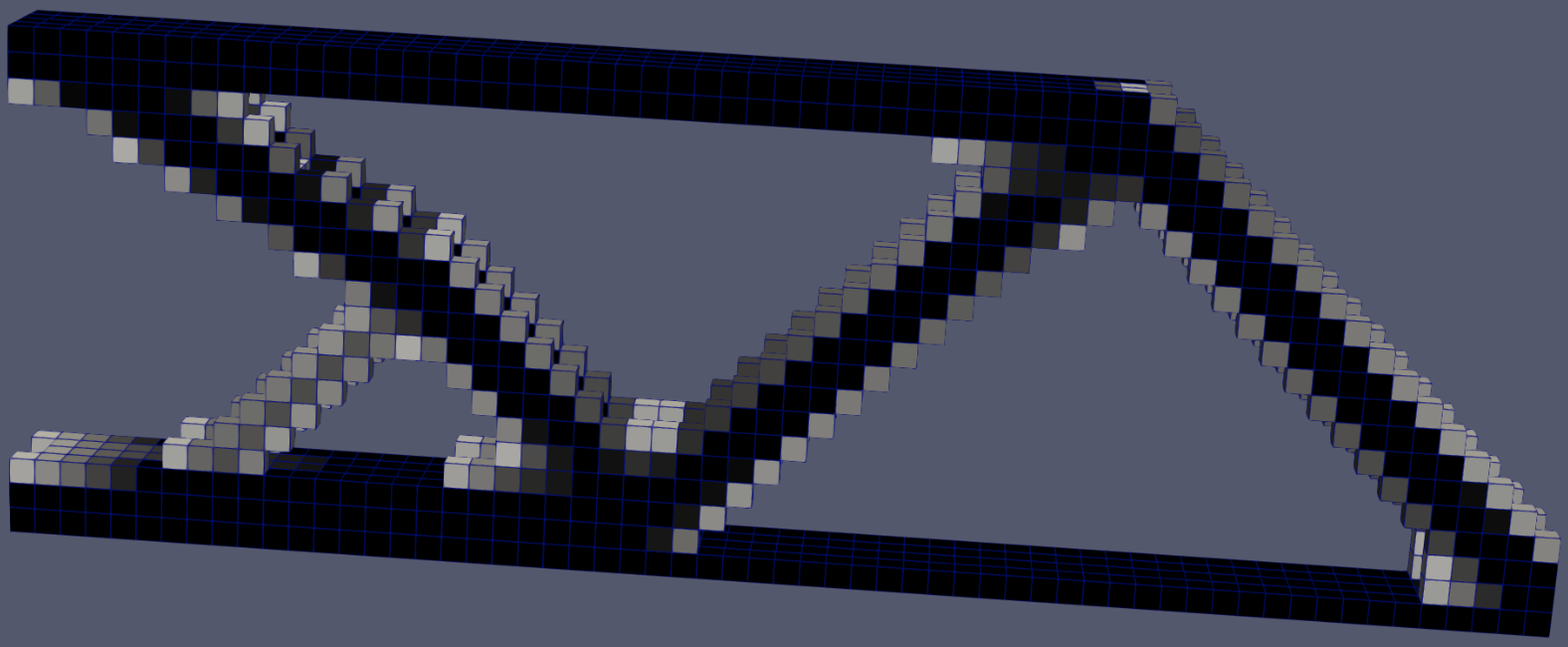}
		\label{fig:canti3d_torch}
	}
	\subfloat[$c(\rho)=2063.5625,~v(\rho)=0.2999$]{
		\includegraphics[width=0.32\textwidth]{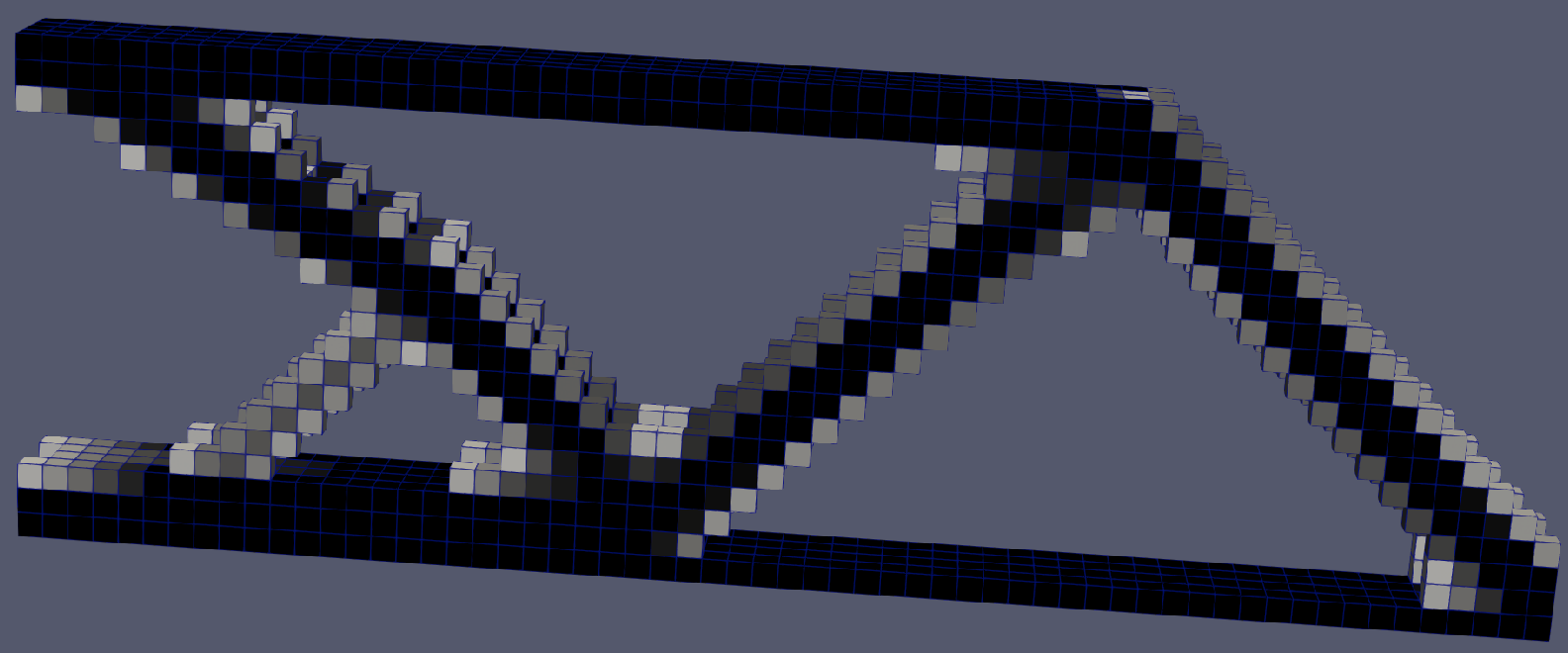}
		\label{fig:canti3d_jax}
	}
	\caption{Final optimized topologies of the 3D cantilever beam using three different backends (elements with $\rho>0.3$ are visualized). Left: NumPy backend; Middle: PyTorch backend; Right: JAX backend.}
	\label{fig:canti3d_backend_compare}
\end{figure}

To validate its consistency and reliability across backends, we conducted tests on the 3D cantilever beam optimization problem described in Sections~\ref{sec:exp_canti3d} and~\ref{sec:exp_canti3d_ad}, using identical parameter settings. The problem was solved independently under each backend, and the results, shown in Figure~\ref{fig:canti3d_backend_compare}, demonstrate high consistency: compliance converges to 2063.5625, the volume fraction stabilizes at 0.3, and the topology layouts are almost identical. This confirms SOPTX's stability and robustness across different computational backends.

In addition to supporting multi-backend switching, SOPTX enables computations on GPU to enhance efficiency, particularly for large-scale 3D TO problems. By default, computations are performed on the CPU (\texttt{device='cpu'}). Users can migrate computations to a GPU (e.g., CUDA devices) in two ways:
\begin{enumerate}
	\item Set the default device globally using the following code:
\begin{lstlisting}[language=python]
bm.set_default_device('cuda')
\end{lstlisting}
	This method is suitable for scenarios where all computations are to be executed on the GPU.
	\item Specify the device during mesh creation, as shown below:
\begin{lstlisting}[language=python]
mesh = UniformMesh(extent=extent, h=h, origin=origin, device='cuda')
\end{lstlisting}
	This approach offers greater flexibility, allowing users to mix CPU and GPU computations within the same program.	
\end{enumerate}

GPU acceleration is essential for 3D TO due to the significant increase in computational complexity as mesh size grows. For instance, in the 3D cantilever beam example from Section~\ref{sec:exp_canti3d}, doubling the mesh to $120\times40\times8$ increases the density design variables to 38,400 and displacement degrees of freedom to 133,947. On the CPU, this leads to prolonged optimization times due to matrix assembly and linear solving costs. In contrast, GPU's parallel computing capabilities greatly enhance efficiency and reduce computation time.

To demonstrate this, we compare CPU and GPU computation times under the PyTorch backend for the same problem. Parameters remain consistent with Section~\ref{sec:exp_canti3d}, except for increasing the filter radius from 1.5 to 3.0 to maintain structural smoothness in the larger domain. Additionally, we use the Conjugate Gradient (CG) method instead of a direct solver (e.g., MUMPS) to handle the larger sparse linear systems, leveraging the GPU's parallel processing for efficient matrix-vector operations.

\begin{figure}[htbp]
	\centering
	\subfloat[$c(\rho)=3498.8266,~v(\rho)=0.3001$]{
		\includegraphics[width=0.5\textwidth]{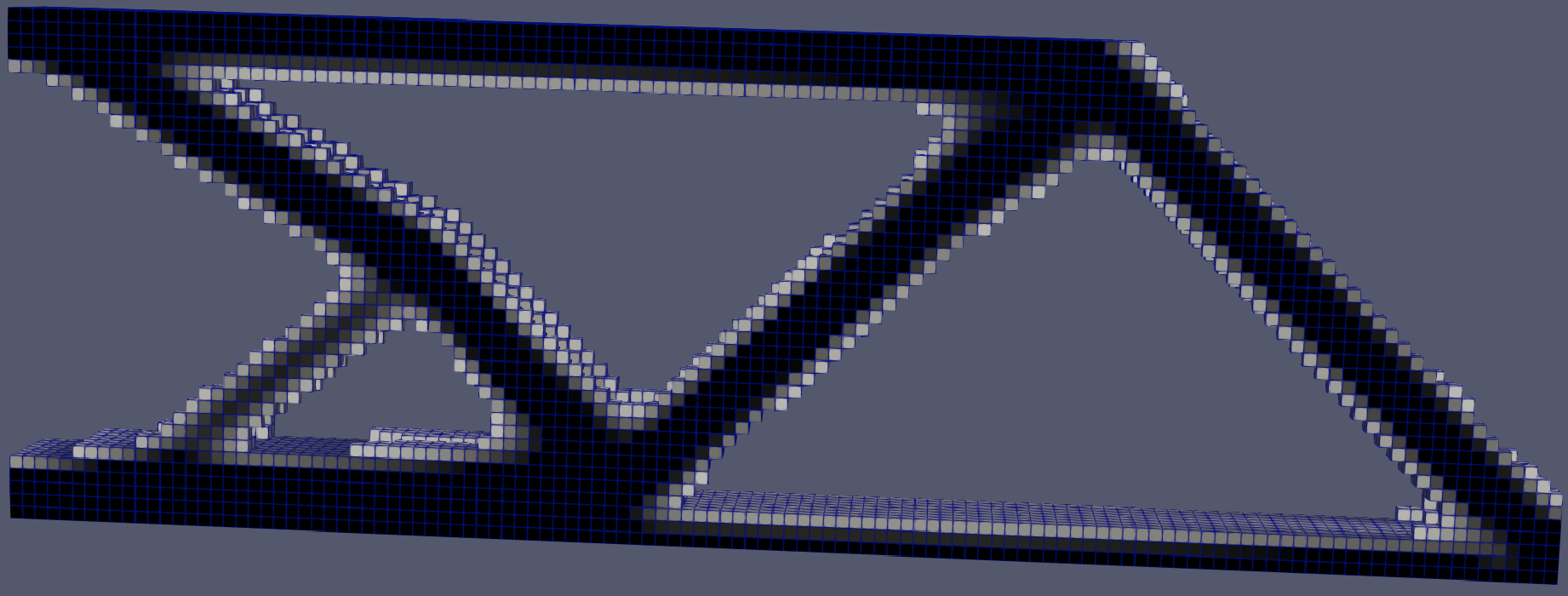}
		\label{fig:canti3d_cpu}
	}
	\subfloat[$c(\rho)=3494.2391,~v(\rho)=0.3000$]{
		\includegraphics[width=0.5\textwidth]{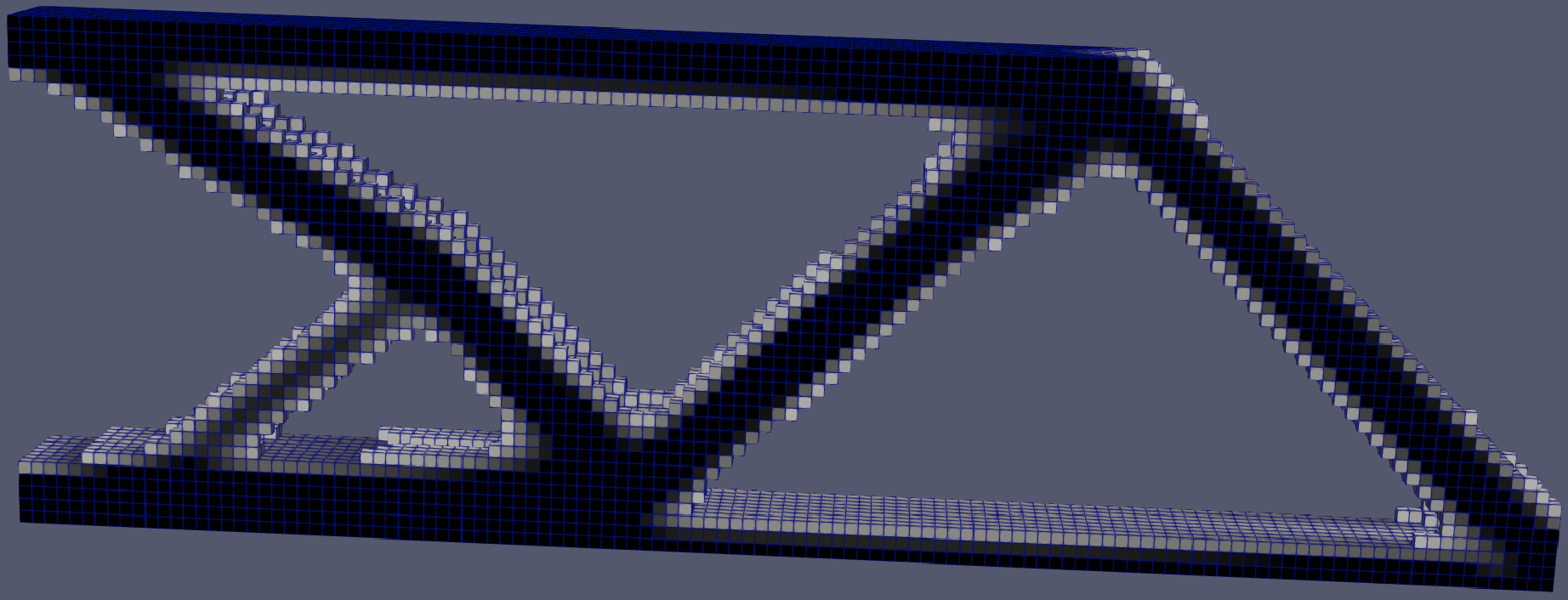}
		\label{fig:canti3d_gpu}
	}
	\caption{Optimized topologies of the 3D cantilever beam using CPU (left) and GPU (right), with elements of density $\rho > 0.3$ visualized.}
	\label{fig:canti3d_device_compare}
\end{figure}

As shown in Figure~\ref{fig:canti3d_device_compare}, the optimized topologies using CPU and GPU are nearly identical. Minor differences in compliance and volume fraction, both within $1\%$, arise from variations in floating-point precision and parallel implementation strategies across hardware platforms. These differences are negligible in TO and do not impact the quality or performance of the final design.

Given the consistency in optimization results, we evaluate the computational efficiency of CPU versus GPU, as detailed in Table~\ref{tab:device_performance}. The GPU acceleration markedly reduces the total optimization time from \SI{3872.041}{s} (CPU) to \SI{479.412}{s} (GPU), achieving an approximately speedup of 8.1 times. This superior performance is evident in both initial and subsequent iterations, highlighting the GPU's parallel computing advantages.
\begin{table}[htbp]
	\centering
	\setlength{\tabcolsep}{4pt} 
	\caption{Performance comparison between CPU and GPU for the 3D cantilever beam optimization problem.}
	\begin{tabular}{ccccc}
		\toprule
		\textbf{Computational Device} & \textbf{Iterations} & \textbf{Total Time (s)} & \textbf{1st Iter. Time (s)} & \textbf{Avg. Iter. Time (s)} \\
		\midrule
		CPU (PyTorch) & 155 & 3872.041 & 15.846 & 25.040 \\
		GPU (PyTorch) & 155 & 479.412 & 2.111 & 3.099 \\
		\bottomrule
	\end{tabular}
	\label{tab:device_performance}
\end{table}

Notably, the Conjugate Gradient (CG) method exhibits faster convergence in early iterations due to the uniform material density, which results in a well-conditioned stiffness matrix. As optimization advances, increased heterogeneity in material distribution elevates the matrix condition number, thereby prolonging CG convergence and iteration times.

\section{Conclusion}
This study presents SOPTX, a high-performance topology optimization (TO) framework built on the open-source FEALPy platform. SOPTX addresses key challenges in TO by offering a modular, efficient, and accessible solution for both academic and industrial applications. Its fully open-source design eliminates reliance on commercial software, broadening its accessibility for research and education.

The core innovations of SOPTX are reflected in the following three aspects:
\begin{enumerate} 
	\item \textbf{Modular Architecture}: A loosely coupled design supports 2D and 3D TO problems across structured and unstructured meshes. This flexibility, paired with a rich set of configurable components, empowers users to tailor and extend optimization workflows efficiently.
	\item \textbf{Multi-Backend Support and Automatic differentiation}: Integration of NumPy, PyTorch, and JAX enables fast central processing units (CPUs) and graphics processing units (GPUs) computation. In benchmark tests, GPU computations reduced time to approximately $12\%$ of CPU time. Automatic differentiation (AD) automates gradient calculations, enhancing accuracy.
	\item \textbf{Efficient Matrix Assembly}: By optimizing assembly processes and exploiting sparsity, SOPTX reduces average assembly time from \SI{0.838}{s} to \SI{0.273}{s} (approximately $33\%$ of the original), significantly boosting efficiency.
\end{enumerate}

The experimental analyses in Sections~\ref{sec:exp_mbb_beam} through~\ref{sec:exp_multi_backend} demonstrate SOPTX's exceptional performance across key areas: seamless model and filter switching, flexible optimization algorithms, 3D problem extension, enhanced computational efficiency, automatic differentiation, and multi-backend support. These results underscore SOPTX's technical strengths and versatility.

Looking forward, SOPTX is set to expand its capabilities:
\begin{itemize} 
	\item \textbf{Level Set Methods}: Integrating FEALPy's level set module will enable boundary-clear optimization, generating smooth, manufacturable designs for precision applications.
	\item \textbf{Adaptive Mesh Refinement}: Incorporating local refinement near boundaries will improve accuracy and detail in complex geometries.
	\item \textbf{Multiphysics Optimization}: Leveraging FEALPy's solvers for heat conduction and Navier Stokes equations, SOPTX can address thermo-fluid-structure coupling, optimizing designs like aerospace thermal protection systems or electronic heat dissipation structures.
	\item \textbf{Manufacturing Constraints}: Future support for stress and additive manufacturing constraints (e.g., overhang control, support minimization) will enhance SOPTX's applicability in aerospace, automotive, and biomedical fields, aligning with lightweight and manufacturability goals.
\end{itemize}
These extensions will solidify SOPTX's role in advancing TO for both research and industry.

\section*{Acknowledgments}
National Natural Science Foundation of China (NSFC) (Grant Nos. 12371410, 12261131501) the construction of innovative provinces in Hunan Province (Grant No. 2021GK1010)

\appendix
\renewcommand{\thesection}{Appendix\,\Alph{section}}
\section{Cantilever2dData1}
\label{sec:code_canti_2d}
\begin{lstlisting}[language=python]
class Cantilever2dData1:
	def __init__(self, 
				xmin: float, xmax: float, 
				ymin: float, ymax: float,
				T: float = -1):
	self.xmin, self.xmax = xmin, xmax
	self.ymin, self.ymax = ymin, ymax
	self.T = T 
	self.eps = 1e-12
	
	def domain(self) -> list:
		box = [self.xmin, self.xmax, self.ymin, self.ymax]
		return box
	
	@cartesian
	def force(self, points: TensorLike) -> TensorLike:
		domain = self.domain()
		x = points[..., 0]
		y = points[..., 1]
		coord = (
			(bm.abs(x - domain[1]) < self.eps) & 
			(bm.abs(y - domain[2]) < self.eps)
			)
		kwargs = bm.context(points)
		val = bm.zeros(points.shape, **kwargs)
		val[coord, 1] = self.T
		return val
	
	@cartesian
	def dirichlet(self, points: TensorLike) -> TensorLike:
		kwargs = bm.context(points)
		return bm.zeros(points.shape, **kwargs)
	
	@cartesian
	def is_dirichlet_boundary_dof_x(self, points: TensorLike) -> TensorLike:
		domain = self.domain()
		x = points[..., 0]
		coord = bm.abs(x - domain[0]) < self.eps
		return coord
	
	@cartesian
	def is_dirichlet_boundary_dof_y(self, points: TensorLike) -> TensorLike:
		domain = self.domain()
		x = points[..., 0]
		coord = bm.abs(x - domain[0]) < self.eps
		return coord    
	
	def threshold(self) -> Tuple[Callable, Callable]:
		return (self.is_dirichlet_boundary_dof_x, 
				self.is_dirichlet_boundary_dof_y)
\end{lstlisting}

\section{MBBBeam2dData1}
\label{sec:code_mbb}
\begin{lstlisting}[language=python]
class MBBBeam2dData1:
def __init__(self, 
			xmin: float=0, xmax: float=60, 
			ymin: float=0, ymax: float=20,
			T: float = -1):
	self.xmin, self.xmax = xmin, xmax
	self.ymin, self.ymax = ymin, ymax
	self.T = T
	self.eps = 1e-12

def domain(self) -> list:
	box = [self.xmin, self.xmax, self.ymin, self.ymax]
	return box

@cartesian
def force(self, points: TensorLike) -> TensorLike:
	domain = self.domain()
	x = points[..., 0]
	y = points[..., 1]
	coord = ((bm.abs(x - domain[0]) < self.eps) & 
			 (bm.abs(y - domain[3]) < self.eps))
	kwargs = bm.context(points)
	val = bm.zeros(points.shape, **kwargs)
	val[coord, 1] = self.T
	return val

@cartesian
def dirichlet(self, points: TensorLike) -> TensorLike:
	kwargs = bm.context(points)
	return bm.zeros(points.shape, **kwargs)

@cartesian
def is_dirichlet_boundary_dof_x(self, points: TensorLike) -> TensorLike:
	domain = self.domain()
	x = points[..., 0]
	coord = bm.abs(x - domain[0]) < self.eps
	return coord

@cartesian
def is_dirichlet_boundary_dof_y(self, points: TensorLike) -> TensorLike:
	domain = self.domain()
	x = points[..., 0]
	y = points[..., 1]
	coord = ((bm.abs(x - domain[1]) < self.eps) &
			 (bm.abs(y - domain[0]) < self.eps))
	return coord

def threshold(self) -> Tuple[Callable, Callable]:
	return (self.is_dirichlet_boundary_dof_x, 
			self.is_dirichlet_boundary_dof_y)
\end{lstlisting}

\section{Cantilever3dData1}
\label{sec:code_canti_3d}
\begin{lstlisting}[language=python]
class Cantilever3dData1:
	def __init__(self,
				xmin: float=0, xmax: float=60, 
				ymin: float=0, ymax: float=20,
				zmin: float=0, zmax: float=4,
				T: float = -1):
		self.xmin, self.xmax = xmin, xmax
		self.ymin, self.ymax = ymin, ymax
		self.zmin, self.zmax = zmin, zmax
		self.T = T 
		self.eps = 1e-12
	
	def domain(self) -> list:
		box = [self.xmin, self.xmax, 
		self.ymin, self.ymax, 
		self.zmin, self.zmax]
		return box
	
	@cartesian
	def force(self, points: TensorLike) -> TensorLike:
		domain = self.domain()
		x = points[..., 0]
		y = points[..., 1]
		z = points[..., 2]
		coord = ((bm.abs(x - domain[1]) < self.eps) & 
				 (bm.abs(y - domain[2]) < self.eps))
		kwargs = bm.context(points)
		val = bm.zeros(points.shape, **kwargs)
		val[coord, 1] = self.T
		return val
	
	@cartesian
	def dirichlet(self, points: TensorLike) -> TensorLike:
		kwargs = bm.context(points)
		return bm.zeros(points.shape, **kwargs)
	
	@cartesian
	def is_dirichlet_boundary_dof_x(self, points: TensorLike) -> TensorLike:
		domain = self.domain()
		x = points[..., 0]
		coord = bm.abs(x - domain[0]) < self.eps
		return coord
	
	@cartesian
	def is_dirichlet_boundary_dof_y(self, points: TensorLike) -> TensorLike:
		domain = self.domain()
		x = points[..., 0]
		coord = bm.abs(x - domain[0]) < self.eps
		return coord
	
	@cartesian
	def is_dirichlet_boundary_dof_z(self, points: TensorLike) -> TensorLike:
		domain = self.domain()
		x = points[..., 0]
		coord = bm.abs(x - domain[0]) < self.eps
		return coord
	
	def threshold(self) -> Tuple[Callable, Callable]:
		return (self.is_dirichlet_boundary_dof_x, 
				self.is_dirichlet_boundary_dof_y,
				self.is_dirichlet_boundary_dof_z)
\end{lstlisting}

\section{Automatic Differentiation}
\label{sec:ad_compliance}
\begin{lstlisting}[language=python]
def _compute_gradient_auto(self, rho: TensorLike, u: Optional[TensorLike] = None) -> TensorLike:
	if u is None:
		u = self._update_u(rho)
	
	ke0 = self.solver.get_base_local_stiffness_matrix()
	cell2dof = self.solver.tensor_space.cell_to_dof()
	ue = u[cell2dof] 
	
	def compliance_contribution(rho_i: float, ue_i: TensorLike, ke0_i: TensorLike) -> float:
		E = self.materials.calculate_elastic_modulus(rho_i)
		c_i = -E * bm.einsum('i, ij, j', ue_i, ke0_i, ue_i)
		return c_i
	
	vmap_grad = bm.vmap(lambda r, u, k:bm.jacrev(lambda x: 
		compliance_contribution(x, u, k))(r))
	dc = vmap_grad(rho, ue, ke0)
	return dc
\end{lstlisting}

\bibliographystyle{abbrv}
\bibliography{paper}


\end{document}